\documentclass[12pt]{article}

\usepackage[utf8]{inputenc}
\usepackage{amssymb}
\usepackage{amsmath}
\usepackage{amsthm}
\usepackage{hyperref}
\usepackage{tikz}
\usepackage{dsfont}
\usepackage{geometry}
\usepackage{enumerate}
\usepackage{float}
\usepackage{caption}
\usepackage{subcaption}
\usepackage{graphicx}
\usepackage{tikz}
\usepackage{tikz-cd}
\usepackage{pgfplots}
\pgfplotsset{
        compat=1.3,
    }
\usepackage{mathtools}
\DeclarePairedDelimiterX{\Bra}[1]{\langle}{\rangle}{#1}
\newcommand{\norm}[1]{{\Vert #1 \Vert }}
\DeclareMathOperator{\spn}{span}

\newcommand{\cT}{\mathcal{T}}

\theoremstyle{theorem}
\newtheorem{theorem}{Theorem}
\newtheorem*{theorem*}{Theorem}
\theoremstyle{proposition}

\newtheorem*{remark*}{Remark}
\theoremstyle{definition}

\newtheorem*{definition*}{Defintion}

\newtheorem{lemma}[theorem]{Lemma}
\theoremstyle{example}

\title{Higher-Order Finite Element Approximation\\ of the Dynamic Laplacian}
\author{Nathanael Schilling\thanks{Center for Mathematics, Technical University of Munich, 85747 Garching, Germany},
    Gary Froyland\thanks{School of Mathematics and Statistics, University of New South Wales, Sydney NSW 2052, Australia},
    Oliver Junge\thanks{Center for Mathematics, Technical University of Munich, 85747 Garching, Germany}}

\begin{document}

\maketitle

\begin{abstract}
The dynamic Laplace operator arises from extending problems of isoperimetry from fixed manifolds to manifolds evolved by general nonlinear dynamics.  Eigenfunctions of this operator are used to identify and track finite-time coherent sets, which physically manifest in fluid flows as jets, vortices, and more complicated structures.
Two robust and efficient finite-element discretisation schemes for numerically computing the dynamic Laplacian were proposed in \cite{froylandjunge18}.
In this work we consider higher-order versions of these two numerical schemes and analyse them experimentally. 
We also prove the numerically computed eigenvalues and eigenvectors converge to the true objects for both schemes under certain assumptions.
We provide an efficient implementation of the higher-order element schemes in an accompanying Julia package.
\end{abstract}

\section{Introduction}

The dynamic Laplacian is a second-order partial differential operator underlying a range of methods for computing finite-time coherent sets in finite-time non-autonomous dynamical systems.
It was introduced in \cite{froyland15} in the context of defining sets that remain coherent in a Lagrangian sense via dynamic isoperimetry.
Coherent sets are time-dependent families of sets whose boundaries remain small relative to the volume of the set as the family evolves according to the nonlinear dynamics;  extensions to weighted, curved manifolds and non-volume-preserving dynamics were made in \cite{FK17}.
Coherent sets are captured by the eigenvectors of the dynamic Laplacian
corresponding to the leading eigenvalues (i.e.\ those closest to $0$).
In particular, level sets of the eigenvector corresponding to the first nontrivial eigenvalue can be used to partition the domain into two
coherent sets.  
A dynamic Cheeger inequality \cite{froyland15,FK17} links this eigenvalue to the ratio of boundary size to volume.
Moreover, if $n$ eigenvalues are close to zero followed by a spectral gap, this forces the eigenvectors to be 
close to linear combinations of indicator functions on an $n$-partition \cite{davies1982pt2}. 
We note that the dynamic Laplacian is the Laplace-Beltrami operator of a weighted manifold 
\cite{karraschkeller} and therefore can be used as the time-independent generator of a diffusion process approximating the given (time-dependent) advection-diffusion process. Eigenvectors corresponding to small eigenvalues decay slowest under this diffusion process and its almost-invariant sets in the sense of \cite{DeJu99}.
Algorithmically, coherent sets can be extracted from the eigenvectors of the dynamic Laplacian via, e.g., clustering techniques \cite{froyland05}, optimising eigenbasis separation \cite{deuflhard2005robust}, optimising sublevel sets \cite{froylandpadberg09}, or sparse eigenbasis approximation \cite{froyland2019sparse}.


In most cases, the dynamic Laplacian eigenproblem must be solved numerically.  To this end,  a scheme based on radial basis functions had been proposed \cite{rbf}, which showed high order of convergence, but suffered from a number of drawbacks like high sensitivity with respect to the radius parameter, a non-real spectrum and non-sparseness of the discretized operator.  In \cite{froylandjunge18}, two finite element schemes were proposed (the ``Cauchy-Green'' (CG) and the ``Transfer Operator'' (TO) approach), which eliminated each of these drawbacks.   

Experimentally, only piecewise linear elements were considered in \cite{froylandjunge18}.
In this paper, we consider higher-order (and in particular quadratic)  elements and analyse convergence properties both theoretically and experimentally.   We provide an efficient implementation in the Julia package \href{https://github.com/CoherentStructures/CoherentStructures.jl}{\texttt{CoherentStructures.jl}}.
We find that using $P^2$ elements can give a higher asymptotic order of convergence compared to $P^1$ elements in the ``CG'' approach.
For the ``CG'' approach, classical theory concerning eigenproblems in FEM applies. We also provide some test cases where using $P^2$ elements can greatly reduce
the amount of information needed to calculate partitions of the domain that show important dynamical features.
The question of convergence in the ``TO'' approaches is more subtle. We
prove convergence of eigenvalues and eigenvectors for a family of TO approaches with $P^1$ elements, but the proof does not give any
insight into the convergence rates that should be expected. For $P^2$ elements
we do not observe asymptotically higher orders of convergence even for the simple example of a one-dimensional shift-map on the torus.
This suggests that using $P^2$ elements does not have substantial benefits when using the ``TO'' approach and that one should use the simpler and well-performing linear $P^1$ elements in the ``TO'' schemes.

\section{The Dynamic Laplacian}

Let $\mathcal{I}\subset\mathbb{R}$ denote a finite subset of time and consider a finite family $(\Omega_t)_{t \in \mathcal I}$ of open bounded subsets of $\mathbb R^d$ with Lipschitz boundary.
For each $t \in \mathcal I$ let  $T_t : \Omega_0 \rightarrow \Omega_t$ be a volume-preserving diffeomorphism.
We assume that $T_t$ is sufficiently regular so that $T_t$ and $T_t^{-1}$ can be smoothly extended to the boundary,
and that $0 \in \mathcal I$ with $T_0$ being the identity. A typical
setting in which these conditions apply are $T_t$ taken to be time-$t$ flow maps of a divergence-free vector-field.

Denote the Laplace operator on $\Omega_t$ by $\Delta_t$ for each $t \in \mathcal I$.
Then the \emph{dynamic Laplacian} (an operator on $\Omega_0$) is given by
\begin{align*}
\Delta^{dyn} := \frac{1}{|\mathcal I|}\sum_{t \in \mathcal I} T_t^* \Delta_t T_{t,*}
\end{align*}
where $T_t^*:L^2(\Omega_t)\to L^2(\Omega_0)$ denotes the pullback by $T_t$ defined by $T_t^*f=f\circ T_t$, and $T_{t,*}:L^2(\Omega_0)\to L^2(\Omega_t)$ is the pushforward defined by $T_{t,*}f = f\circ T_t^{-1}$.  

Standard PDE-theoretic arguments can be used to show that $\Delta^{dyn}$ is a uniformly elliptic second-order partial-differential operator \cite{evans,froyland15}, with weak form  \cite{froylandjunge18}
\begin{align}\label{eq:weak_form_CG}
a(u,v) &= \frac{1}{|\mathcal{I}|}\sum_{t \in \mathcal I}\int_\Omega \nabla u \cdot [DT_t]^{-1}([DT_t]^{-1})^T \nabla v \; d\ell^d,
\end{align}
where $\ell^d$ is the $d$-dimensional Lebesgue measure.
The time set $\mathcal I$ can also be a compact interval, and the dynamic Laplacian can be defined by means of an integral over $\mathcal I$ \cite{froyland15,FK17}
$$\Delta^{dyn} := \frac{1}{|\mathcal I|}\int_\mathcal I T_t^* \Delta_t T_{t,*}\ dt,$$ where $|\mathcal I|$ now denotes the length of the interval $\mathcal I$, we do not consider this generalization further here.

The uniform ellipticity of the dynamic Laplacian ensures that the bilinear form in (\ref{eq:weak_form_CG})
is coercive under suitable boundary conditions. These are determined by the choice of the underlying space $S$ on which (\ref{eq:weak_form_CG}) acts.
For natural (Neumann) boundary conditions, we look at $a$ on $\hat{H^1} \times \hat{H^1}$, where $\hat{H}^1 = \hat{H}^1(\Omega_0) := \{v \in H^1(\Omega_0) : \int f \; d\ell^d = 0\}$ is the space of mean-free $H^1$ functions.
For homogeneous Dirichlet boundary conditions, $a$ must be taken to act on $H^1_0 \times H^1_0$.
With either choice of boundary conditions, the variational form of the of the eigenproblem $\Delta^{dyn} u = \lambda u$
becomes
\begin{align}\label{eq:weak_eigenproblem}
a(u,v) = \lambda\Bra{u,v}_{L^2}\quad \mbox{for all}\ v \in S.
\end{align}
In both cases, ellipticity ensures that there exists a countable sequence of pairwise orthogonal eigenvectors $u_0,u_1,\ldots$ corresponding to real eigenvalues $0\geq\lambda_0\geq\lambda_1\geq\cdots$. Furthermore, the span of the eigenvectors is dense in $L^2$ (Dirichlet) and mean free $L^2$ functions (Neumann) \cite{evans,froyland15}.

\subsection{Discretisation with finite elements}

A natural discretisation of the eigenproblem \eqref{eq:weak_eigenproblem} is by using a finite element method (FEM) \cite{froylandjunge18}.  In the standard Galerkin discretization of \eqref{eq:weak_eigenproblem}, $S' \subset S$ is taken to be a finite dimensional approximation space spanned by some basis $(\varphi_1, \dots, \varphi_n)$.  We now find approximate solutions $u\in S'$, $\lambda\in\mathbb{R}$, that satisfy (\ref{eq:weak_eigenproblem}) with $S'$ taking the place of $S$.  In matrix form, the coefficients $\textbf{u} = (\textbf{u}_1,\ldots,\textbf{u}_n)$ of $u$ with respect to the basis $(\varphi_1,\dots,\varphi_n)$  are found by solving the generalized eigenvalue problem
\begin{align}\label{eq:matrixform}
D\textbf{u} = \lambda M \textbf{u}.
\end{align}
Here, $D=(D_{i,j}) = a(\varphi_i,\varphi_j)$ is referred to as the \emph{stiffness matrix} and $M=(M_{i,j}) = \Bra{\varphi_i,\varphi_j}_{L^2}$
as the \emph{mass matrix}.  We use $P^k$ Lagrange nodal basis functions on some (triangular or quadrilateral) mesh for $\varphi_1, \dots, \varphi_n$, but other choices are of also possible in general (see \cite{ernguermond}).

Using a finite element method for approximating $\Delta^{dyn}$ has a number of advantages: The matrix formulation \eqref{eq:matrixform} inherits self-adjointness from the continuous problem and thus the computed eigenpairs are always real. With a suitably localized basis, $D$ and $M$ are both sparse. Sparse Hermitian generalized eigenproblems are well known in the literature, and a number of algorithms exist for efficiently solving them \cite{bai2000templates}.  Also, finite element methods have been widely studied, and there are a range of theoretical results regarding convergence (cf \cite{ernguermond}) that are applicable to some solution approaches.

While the entries of the mass matrix can be computed exactly, the integrand in (\ref{eq:weak_form_CG}) in the entries of the stiffness matrix can be of extremely high variation locally, making an accurate computation possibly expensive.  Several approaches for approximating $D$ were suggested in \cite{froylandjunge18}, which we investigate further in the following sections.

 \section{Convergence for the CG approach}
\label{sec:CGconv}
In the ``CG approach'' of \cite{froylandjunge18}, a quadrature formula is used to approximate (\ref{eq:weak_form_CG}). This is the standard
way of solving elliptic eigenproblems with finite elements.

We work on $P^k$-Lagrange finite element spaces for a family of quasi-uniform meshes
\footnote{A family $\{\mathcal T_h\}_{h>0}$ of meshes is \emph{quasi-uniform}, if it is shape-regular and there exists a constant $c>0$ so that for any $h > 0$ and any element $K \in \mathcal T_h$ one has that $h_K \geq ch$, where $h_K$ is the diameter of the element $K$. Shape regularity
means that there is a constant $\sigma_0$ so that uniformly for all $K$ it holds that $\sigma = \frac{h_K}{\rho_K} \leq \sigma_0$. Here $\rho_K$ is the diameter of the largest ball that can be inscribed in $K$\cite{ernguermond}.
}
triangulations (in $2$ or $3$ dimensions) $\cT_h^0$ of $\Omega_0$, where the family parameter $h$ is related to the mesh size.  Let $\varphi^{1}_h,\ldots,\varphi^{N(h)}_h$ be the associated nodal basis, $S_h^0 := \mbox{span}\{\varphi^{1}_h,\ldots,\varphi^{N(h)}_h\}$ the finite element approximation space,
and $\hat S_h^0 := \hat H^1 \cap S_h^0$.  Moreover, let $ \lambda_{i,h}$ be the  Ritz-values with corresponding Ritz-vectors $v_{i,h}$ in $\hat S_h^0$, i.e. $\lambda_{i,h}$ is the minimum value of $a(v,v)$ over the set $\{ v \in \hat S_h^0 ;  \norm{v}_{L^2} = 1 \Bra{v, v_j} = 0 \mbox{ for } j = 1,\cdots,(i-1) \}$ and $v_{i,h}$ is the corresponding minimizer.

Under the assumption that the eigenvectors of $a:\hat H^1\times \hat H^1\to\mathbb{R}$ are in $H^{k+1}$, classical results from FEM theory \cite{strangfix,ernguermond} give that (where $C$ stands for a generic constant):
\begin{align}
\label{eq:CG_lambda_error}
\frac{\lambda_{i,h} - \lambda_i}{\lambda_i} \leq C h^{2k} \qquad (h\to 0).
\end{align}
Moreoever, it holds for simple (and $L^2$-normalised) eigenvectors that
\begin{align}
\label{eq:CG_v_error_1}
\norm{v_i - v_{i,h}}_{H^1} \leq C h^{k} \qquad (h\to 0).
\end{align}
In the case that the problem (\ref{eq:weak_form_CG}) satisfies the technical assumption of being \emph{regularizing} \cite{ernguermond} (i.e $\norm{Pu}_{H^2} \leq C \norm{u}_{L^2}$), the ``Aubin-Nitsche trick'' gives
\begin{align}
\label{eq:CG_v_error_2}
\norm{v_i - v_{i,h}}_{L^2} \leq C h^{k+1} \qquad (h\to 0).
\end{align}
Note that not all elliptic problems are regularizing. We thus cannot expect \eqref{eq:CG_v_error_2} to hold in general.

If the entries of the stiffness and the mass matrix are approximated with a quadrature rule of order $2k-1$, the convergence orders are unaffected, provided that the quadrature points contain a $\mathcal P^k$ unisolvent set \cite{quadraturepaper}. Moreoever, the results generalize to non-simple eigenspaces and Kato's subspace distance used for the error\cite{quadraturepaper}.

\subsection{Numerical Experiments}

We now aim at reproducing the predicted convergence rates in numerical experiments.  As a reference, we compute ($L^2$ normalized) eigenpairs $(\lambda_i,v_i)$ on a very fine mesh and estimate the error in the eigenvector $v_{i,h}$ by computing the $L^2$-distance of $v_{i,h}$ to the closest reference eigenvector, given by the expression
\begin{align}\label{eq:error_formula}
e_{i,h} := \sqrt{1 - |\Bra{v_i,v_{i,h}}_{L^2}|^2}.
\end{align}
The $L^2$ inner-product is approximated by first interpolating to the fine grid,
and then calculating the inner product there using quadrature.

To generalize this to higher-dimensional eigenspaces, let $\tilde V$ and $V$ be two $m$-dimensional subspaces of $L^2$ with orthonormal bases given by $\{\tilde v_1, \dots , \tilde v_m\}$ and $\{v_1,\dots,v_m\}$. Let $\tilde P, P$ be orthogonal projections onto $\tilde V$ and $V$ respectively.
As a measure of the subspace error \cite{kato}, we maximize $\norm{\tilde v - P\tilde v}_{L^2}$ over
$\tilde v \in \tilde V$ with $\norm{\tilde v}_{L^2} = 1$.
This is equivalent to maximizing $\sqrt{1 - \norm{P\tilde v}^2}$,
which has maximum given by $\sqrt{1 - \rho_{min}(P)^2}$ where $\rho_{min}$ is the smallest singular value of $P$ on $\tilde V$ 
which is also that of the matrix $(\Bra{\tilde v_i,v_j}_{L^2})_{i,j=1}^m$.

In general, one is not directly interested in the eigenvectors themselves, but in a partition of the domain obtained by a suitable post-processing of the eigenvectors (e.g., as mentioned, thresholding, clustering or sparse eigenbases).  This motivates us making a qualitative comparison of such partitions based on eigenvectors for different approximation spaces.  Here, we focus in particular on the question of how well Lagrangian coherent sets can be found with as little information about the flow as possible. Evidently, in this low data case some features may not be accurately resolved even if the eigenvectors clearly give some indication of their existence.
We therefore also look at how coarse the grid can be made without affecting the topology of such a partition in some test cases.

In the sequel, mesh sizes of the form $n\times m$ refer to a regular triangular mesh with $n \times m$ (non-inner) nodes in each direction. 

For all time-integration done in this work we used a relative and absolute tolerance of $10^{-8}$ and the \texttt{DifferentialEquations.jl} Julia package \cite{juliadiffeq} with the \texttt{BS5()} solver \cite{Bogacki1996}.  Derivatives of flow maps were approximated with second order central finite differences.
The stiffness matrix was calculated by approximating the integral in \eqref{eq:weak_form_CG} on each element with nodal basis functions for $u$ and $v$ using quadrature.
The results from this were then summed over all elements in the support of a basis function to give the corresponding entry in the stiffness matrix. 
Similarly, the integral in the mass matrix was also calculated element-wise and then additively combined to give the full mass matrix. 
The \texttt{CoherentStructures.jl} package internally uses the \texttt{JuAFEM.jl} package \cite{juafem}.

\subsubsection{Standard Map}

\begin{figure}[H]
\begin{subfigure}{0.5\linewidth}
\begin{tikzpicture}[]
\begin{axis}[scale only axis,width=0.8\linewidth, legend pos = {north west}, ylabel = {Relative eigenvalue error}, xmin = {0.015641212122425755}, xmax = {0.6162109095189517}, ymax = {0.10859816912680587}, ymode = {log}, xlabel = {Mesh width}, unbounded coords=jump,scaled x ticks = false,xlabel style = {font = {\fontsize{11 pt}{14.3 pt}\selectfont}, color = {rgb,1:red,0.00000000;green,0.00000000;blue,0.00000000}, draw opacity = 1.0, rotate = 0.0},log basis x=10,xmajorgrids = true,xtick = {0.015848931924611134,0.039810717055349734,0.1,0.251188643150958},xticklabels = {$10^{-1.8}$,$10^{-1.4}$,$10^{-1.0}$,$10^{-0.6}$},xtick align = inside,xticklabel style = {font = {\fontsize{8 pt}{10.4 pt}\selectfont}, color = {rgb,1:red,0.00000000;green,0.00000000;blue,0.00000000}, draw opacity = 1.0, rotate = 0.0},x grid style = {color = {rgb,1:red,0.00000000;green,0.00000000;blue,0.00000000},
draw opacity = 0.1,
line width = 0.5,
solid},axis lines* = left,x axis line style = {color = {rgb,1:red,0.00000000;green,0.00000000;blue,0.00000000},
draw opacity = 1.0,
line width = 1,
solid},scaled y ticks = false,ylabel style = {font = {\fontsize{11 pt}{14.3 pt}\selectfont}, color = {rgb,1:red,0.00000000;green,0.00000000;blue,0.00000000}, draw opacity = 1.0, rotate = 0.0},log basis y=10,ymajorgrids = true,ytick = {1.0e-8,1.0e-6,9.999999999999999e-5,0.01},yticklabels = {$10^{-8}$,$10^{-6}$,$10^{-4}$,$10^{-2}$},ytick align = inside,yticklabel style = {font = {\fontsize{8 pt}{10.4 pt}\selectfont}, color = {rgb,1:red,0.00000000;green,0.00000000;blue,0.00000000}, draw opacity = 1.0, rotate = 0.0},y grid style = {color = {rgb,1:red,0.00000000;green,0.00000000;blue,0.00000000},
draw opacity = 0.1,
line width = 0.5,
solid},axis lines* = left,y axis line style = {color = {rgb,1:red,0.00000000;green,0.00000000;blue,0.00000000},
draw opacity = 1.0,
line width = 1,
solid},    xshift = 0.0mm,
    yshift = 0.0mm,
    axis background/.style={fill={rgb,1:red,1.00000000;green,1.00000000;blue,1.00000000}}
,legend style = {color = {rgb,1:red,0.00000000;green,0.00000000;blue,0.00000000},
draw opacity = 1.0,
line width = 1,
solid,fill = {rgb,1:red,1.00000000;green,1.00000000;blue,1.00000000},font = {\fontsize{8 pt}{10.4 pt}\selectfont}},colorbar style={title=}, xmode = {log}, ymin = {1.9897937439242972e-10} ]\addplot+[draw=none, color = {rgb,1:red,0.00000000;green,0.50196078;blue,0.00000000},
draw opacity = 1.0,
line width = 0,
solid,mark = diamond*,
mark size = 3.0,
mark options = {
    color = {rgb,1:red,0.00000000;green,0.00000000;blue,0.00000000}, draw opacity = 1.0,
    fill = {rgb,1:red,0.00000000;green,0.50196078;blue,0.00000000}, fill opacity = 1.0,
    line width = 1,
    rotate = 0,
    solid
}] coordinates {
(0.5553603672697958, 0.0004018639199799052)
};
\addlegendentry{P2 elements (4.0)}
\addlegendentry{P2 elements (4.0)}
\addlegendentry{P2 elements (4.0)}
\addlegendentry{P2 elements (4.0)}
\addlegendentry{P2 elements (4.0)}
\addlegendentry{P2 elements (4.0)}
\addplot+[draw=none, color = {rgb,1:red,0.00000000;green,0.50196078;blue,0.00000000},
draw opacity = 1.0,
line width = 0,
solid,mark = diamond*,
mark size = 3.0,
mark options = {
    color = {rgb,1:red,0.00000000;green,0.00000000;blue,0.00000000}, draw opacity = 1.0,
    fill = {rgb,1:red,0.00000000;green,0.50196078;blue,0.00000000}, fill opacity = 1.0,
    line width = 1,
    rotate = 0,
    solid
},forget plot] coordinates {
(0.2776801836348979, 2.717911043324043e-5)
};
\addplot+[draw=none, color = {rgb,1:red,0.00000000;green,0.50196078;blue,0.00000000},
draw opacity = 1.0,
line width = 0,
solid,mark = diamond*,
mark size = 3.0,
mark options = {
    color = {rgb,1:red,0.00000000;green,0.00000000;blue,0.00000000}, draw opacity = 1.0,
    fill = {rgb,1:red,0.00000000;green,0.50196078;blue,0.00000000}, fill opacity = 1.0,
    line width = 1,
    rotate = 0,
    solid
},forget plot] coordinates {
(0.13884009181744894, 1.740317427634017e-6)
};
\addplot+[draw=none, color = {rgb,1:red,0.00000000;green,0.50196078;blue,0.00000000},
draw opacity = 1.0,
line width = 0,
solid,mark = diamond*,
mark size = 3.0,
mark options = {
    color = {rgb,1:red,0.00000000;green,0.00000000;blue,0.00000000}, draw opacity = 1.0,
    fill = {rgb,1:red,0.00000000;green,0.50196078;blue,0.00000000}, fill opacity = 1.0,
    line width = 1,
    rotate = 0,
    solid
},forget plot] coordinates {
(0.06942004590872447, 1.094107666364215e-7)
};
\addplot+[draw=none, color = {rgb,1:red,0.00000000;green,0.50196078;blue,0.00000000},
draw opacity = 1.0,
line width = 0,
solid,mark = diamond*,
mark size = 3.0,
mark options = {
    color = {rgb,1:red,0.00000000;green,0.00000000;blue,0.00000000}, draw opacity = 1.0,
    fill = {rgb,1:red,0.00000000;green,0.50196078;blue,0.00000000}, fill opacity = 1.0,
    line width = 1,
    rotate = 0,
    solid
},forget plot] coordinates {
(0.034710022954362235, 6.778806058655367e-9)
};
\addplot+[draw=none, color = {rgb,1:red,0.00000000;green,0.50196078;blue,0.00000000},
draw opacity = 1.0,
line width = 0,
solid,mark = diamond*,
mark size = 3.0,
mark options = {
    color = {rgb,1:red,0.00000000;green,0.00000000;blue,0.00000000}, draw opacity = 1.0,
    fill = {rgb,1:red,0.00000000;green,0.50196078;blue,0.00000000}, fill opacity = 1.0,
    line width = 1,
    rotate = 0,
    solid
},forget plot] coordinates {
(0.017355011477181118, 3.516268699270139e-10)
};
\addplot+[draw=none, color = {rgb,1:red,0.00000000;green,0.50196078;blue,0.00000000},
draw opacity = 1.0,
line width = 0,
solid,mark = triangle*,
mark size = 3.0,
mark options = {
    color = {rgb,1:red,0.00000000;green,0.00000000;blue,0.00000000}, draw opacity = 1.0,
    fill = {rgb,1:red,0.00000000;green,0.50196078;blue,0.00000000}, fill opacity = 1.0,
    line width = 1,
    rotate = 0,
    solid
}] coordinates {
(0.5553603672697958, 0.05924545441242054)
};
\addlegendentry{P1 elements (2.0)}
\addlegendentry{P1 elements (2.0)}
\addlegendentry{P1 elements (2.0)}
\addlegendentry{P1 elements (2.0)}
\addlegendentry{P1 elements (2.0)}
\addlegendentry{P1 elements (2.0)}
\addplot+[draw=none, color = {rgb,1:red,0.00000000;green,0.50196078;blue,0.00000000},
draw opacity = 1.0,
line width = 0,
solid,mark = triangle*,
mark size = 3.0,
mark options = {
    color = {rgb,1:red,0.00000000;green,0.00000000;blue,0.00000000}, draw opacity = 1.0,
    fill = {rgb,1:red,0.00000000;green,0.50196078;blue,0.00000000}, fill opacity = 1.0,
    line width = 1,
    rotate = 0,
    solid
},forget plot] coordinates {
(0.2776801836348979, 0.015203604963563743)
};
\addplot+[draw=none, color = {rgb,1:red,0.00000000;green,0.50196078;blue,0.00000000},
draw opacity = 1.0,
line width = 0,
solid,mark = triangle*,
mark size = 3.0,
mark options = {
    color = {rgb,1:red,0.00000000;green,0.00000000;blue,0.00000000}, draw opacity = 1.0,
    fill = {rgb,1:red,0.00000000;green,0.50196078;blue,0.00000000}, fill opacity = 1.0,
    line width = 1,
    rotate = 0,
    solid
},forget plot] coordinates {
(0.13884009181744894, 0.003830485900351286)
};
\addplot+[draw=none, color = {rgb,1:red,0.00000000;green,0.50196078;blue,0.00000000},
draw opacity = 1.0,
line width = 0,
solid,mark = triangle*,
mark size = 3.0,
mark options = {
    color = {rgb,1:red,0.00000000;green,0.00000000;blue,0.00000000}, draw opacity = 1.0,
    fill = {rgb,1:red,0.00000000;green,0.50196078;blue,0.00000000}, fill opacity = 1.0,
    line width = 1,
    rotate = 0,
    solid
},forget plot] coordinates {
(0.06942004590872447, 0.0009595896941673767)
};
\addplot+[draw=none, color = {rgb,1:red,0.00000000;green,0.50196078;blue,0.00000000},
draw opacity = 1.0,
line width = 0,
solid,mark = triangle*,
mark size = 3.0,
mark options = {
    color = {rgb,1:red,0.00000000;green,0.00000000;blue,0.00000000}, draw opacity = 1.0,
    fill = {rgb,1:red,0.00000000;green,0.50196078;blue,0.00000000}, fill opacity = 1.0,
    line width = 1,
    rotate = 0,
    solid
},forget plot] coordinates {
(0.034710022954362235, 0.0002400225654180785)
};
\addplot+[draw=none, color = {rgb,1:red,0.00000000;green,0.50196078;blue,0.00000000},
draw opacity = 1.0,
line width = 0,
solid,mark = triangle*,
mark size = 3.0,
mark options = {
    color = {rgb,1:red,0.00000000;green,0.00000000;blue,0.00000000}, draw opacity = 1.0,
    fill = {rgb,1:red,0.00000000;green,0.50196078;blue,0.00000000}, fill opacity = 1.0,
    line width = 1,
    rotate = 0,
    solid
},forget plot] coordinates {
(0.017355011477181118, 6.0013444384765215e-5)
};
\addplot+ [color = {rgb,1:red,0.00000000;green,0.50196078;blue,0.00000000},
draw opacity = 1.0,
line width = 1,
solid,mark = none,
mark size = 2.0,
mark options = {
    color = {rgb,1:red,0.00000000;green,0.00000000;blue,0.00000000}, draw opacity = 1.0,
    fill = {rgb,1:red,0.00000000;green,0.50196078;blue,0.00000000}, fill opacity = 1.0,
    line width = 1,
    rotate = 0,
    solid
},forget plot]coordinates {
(0.5553603672697958, 0.061453767049999476)
(0.2776801836348979, 0.015363441762499874)
(0.13884009181744894, 0.0038408604406249694)
(0.06942004590872447, 0.000960215110156242)
(0.034710022954362235, 0.00024005377753906045)
(0.017355011477181118, 6.001344438476516e-5)
};
\addplot+ [color = {rgb,1:red,0.00000000;green,0.50196078;blue,0.00000000},
draw opacity = 1.0,
line width = 1,
solid,mark = none,
mark size = 2.0,
mark options = {
    color = {rgb,1:red,0.00000000;green,0.00000000;blue,0.00000000}, draw opacity = 1.0,
    fill = {rgb,1:red,0.00000000;green,0.50196078;blue,0.00000000}, fill opacity = 1.0,
    line width = 1,
    rotate = 0,
    solid
},forget plot]coordinates {
(0.5553603672697958, 0.0003687074967605875)
(0.2776801836348979, 2.3044218547536705e-5)
(0.13884009181744894, 1.4402636592210464e-6)
(0.06942004590872447, 9.001647870131535e-8)
(0.034710022954362235, 5.626029918832217e-9)
(0.017355011477181118, 3.5162686992701343e-10)
};
\end{axis}

\end{tikzpicture}
\end{subfigure}
\begin{subfigure}{0.5\linewidth}
\begin{tikzpicture}[]
\begin{axis}[scale only axis,width=0.8\linewidth, legend pos = {north west}, ylabel = {Eigenspace  error}, xmin = {0.015641212122425755}, xmax = {0.6162109095189517}, ymax = {0.048919113372689664}, ymode = {log}, xlabel = {Mesh width}, unbounded coords=jump,scaled x ticks = false,xlabel style = {font = {\fontsize{11 pt}{14.3 pt}\selectfont}, color = {rgb,1:red,0.00000000;green,0.00000000;blue,0.00000000}, draw opacity = 1.0, rotate = 0.0},log basis x=10,xmajorgrids = true,xtick = {0.015848931924611134,0.039810717055349734,0.1,0.251188643150958},xticklabels = {$10^{-1.8}$,$10^{-1.4}$,$10^{-1.0}$,$10^{-0.6}$},xtick align = inside,xticklabel style = {font = {\fontsize{8 pt}{10.4 pt}\selectfont}, color = {rgb,1:red,0.00000000;green,0.00000000;blue,0.00000000}, draw opacity = 1.0, rotate = 0.0},x grid style = {color = {rgb,1:red,0.00000000;green,0.00000000;blue,0.00000000},
draw opacity = 0.1,
line width = 0.5,
solid},axis lines* = left,x axis line style = {color = {rgb,1:red,0.00000000;green,0.00000000;blue,0.00000000},
draw opacity = 1.0,
line width = 1,
solid},scaled y ticks = false,ylabel style = {font = {\fontsize{11 pt}{14.3 pt}\selectfont}, color = {rgb,1:red,0.00000000;green,0.00000000;blue,0.00000000}, draw opacity = 1.0, rotate = 0.0},log basis y=10,ymajorgrids = true,ytick = {1.0e-7,1.0e-6,1.0e-5,9.999999999999999e-5,0.001,0.01},yticklabels = {$10^{-7}$,$10^{-6}$,$10^{-5}$,$10^{-4}$,$10^{-3}$,$10^{-2}$},ytick align = inside,yticklabel style = {font = {\fontsize{8 pt}{10.4 pt}\selectfont}, color = {rgb,1:red,0.00000000;green,0.00000000;blue,0.00000000}, draw opacity = 1.0, rotate = 0.0},y grid style = {color = {rgb,1:red,0.00000000;green,0.00000000;blue,0.00000000},
draw opacity = 0.1,
line width = 0.5,
solid},axis lines* = left,y axis line style = {color = {rgb,1:red,0.00000000;green,0.00000000;blue,0.00000000},
draw opacity = 1.0,
line width = 1,
solid},    xshift = 0.0mm,
    yshift = 0.0mm,
    axis background/.style={fill={rgb,1:red,1.00000000;green,1.00000000;blue,1.00000000}}
,legend style = {color = {rgb,1:red,0.00000000;green,0.00000000;blue,0.00000000},
draw opacity = 1.0,
line width = 1,
solid,fill = {rgb,1:red,1.00000000;green,1.00000000;blue,1.00000000},font = {\fontsize{8 pt}{10.4 pt}\selectfont}},colorbar style={title=}, xmode = {log}, ymin = {3.1476706435933934e-8}]\addplot+[draw=none, color = {rgb,1:red,0.00000000;green,0.50196078;blue,0.00000000},
draw opacity = 1.0,
line width = 0,
solid,mark = diamond*,
mark size = 3.0,
mark options = {
    color = {rgb,1:red,0.00000000;green,0.00000000;blue,0.00000000}, draw opacity = 1.0,
    fill = {rgb,1:red,0.00000000;green,0.50196078;blue,0.00000000}, fill opacity = 1.0,
    line width = 1,
    rotate = 0,
    solid
}] coordinates {
(0.5553603672697958, 0.0009534084281314785)
};
\addlegendentry{P2 elements (3.0)}
\addlegendentry{P2 elements (3.0)}
\addlegendentry{P2 elements (3.0)}
\addlegendentry{P2 elements (3.0)}
\addlegendentry{P2 elements (3.0)}
\addlegendentry{P2 elements (3.0)}
\addplot+[draw=none, color = {rgb,1:red,0.00000000;green,0.50196078;blue,0.00000000},
draw opacity = 1.0,
line width = 0,
solid,mark = diamond*,
mark size = 3.0,
mark options = {
    color = {rgb,1:red,0.00000000;green,0.00000000;blue,0.00000000}, draw opacity = 1.0,
    fill = {rgb,1:red,0.00000000;green,0.50196078;blue,0.00000000}, fill opacity = 1.0,
    line width = 1,
    rotate = 0,
    solid
},forget plot] coordinates {
(0.2776801836348979, 0.00011321635353271144)
};
\addplot+[draw=none, color = {rgb,1:red,0.00000000;green,0.50196078;blue,0.00000000},
draw opacity = 1.0,
line width = 0,
solid,mark = diamond*,
mark size = 3.0,
mark options = {
    color = {rgb,1:red,0.00000000;green,0.00000000;blue,0.00000000}, draw opacity = 1.0,
    fill = {rgb,1:red,0.00000000;green,0.50196078;blue,0.00000000}, fill opacity = 1.0,
    line width = 1,
    rotate = 0,
    solid
},forget plot] coordinates {
(0.13884009181744894, 1.3812950443028499e-5)
};
\addplot+[draw=none, color = {rgb,1:red,0.00000000;green,0.50196078;blue,0.00000000},
draw opacity = 1.0,
line width = 0,
solid,mark = diamond*,
mark size = 3.0,
mark options = {
    color = {rgb,1:red,0.00000000;green,0.00000000;blue,0.00000000}, draw opacity = 1.0,
    fill = {rgb,1:red,0.00000000;green,0.50196078;blue,0.00000000}, fill opacity = 1.0,
    line width = 1,
    rotate = 0,
    solid
},forget plot] coordinates {
(0.06942004590872447, 1.7138278892433494e-6)
};
\addplot+[draw=none, color = {rgb,1:red,0.00000000;green,0.50196078;blue,0.00000000},
draw opacity = 1.0,
line width = 0,
solid,mark = diamond*,
mark size = 3.0,
mark options = {
    color = {rgb,1:red,0.00000000;green,0.00000000;blue,0.00000000}, draw opacity = 1.0,
    fill = {rgb,1:red,0.00000000;green,0.50196078;blue,0.00000000}, fill opacity = 1.0,
    line width = 1,
    rotate = 0,
    solid
},forget plot] coordinates {
(0.034710022954362235, 2.1773039096895542e-7)
};
\addplot+[draw=none, color = {rgb,1:red,0.00000000;green,0.50196078;blue,0.00000000},
draw opacity = 1.0,
line width = 0,
solid,mark = diamond*,
mark size = 3.0,
mark options = {
    color = {rgb,1:red,0.00000000;green,0.00000000;blue,0.00000000}, draw opacity = 1.0,
    fill = {rgb,1:red,0.00000000;green,0.50196078;blue,0.00000000}, fill opacity = 1.0,
    line width = 1,
    rotate = 0,
    solid
},forget plot] coordinates {
(0.017355011477181118, 4.712160915387242e-8)
};
\addplot+[draw=none, color = {rgb,1:red,0.00000000;green,0.50196078;blue,0.00000000},
draw opacity = 1.0,
line width = 0,
solid,mark = triangle*,
mark size = 3.0,
mark options = {
    color = {rgb,1:red,0.00000000;green,0.00000000;blue,0.00000000}, draw opacity = 1.0,
    fill = {rgb,1:red,0.00000000;green,0.50196078;blue,0.00000000}, fill opacity = 1.0,
    line width = 1,
    rotate = 0,
    solid
}] coordinates {
(0.5553603672697958, 0.029290200782601434)
};
\addlegendentry{P1 elements (2.0)}
\addlegendentry{P1 elements (2.0)}
\addlegendentry{P1 elements (2.0)}
\addlegendentry{P1 elements (2.0)}
\addlegendentry{P1 elements (2.0)}
\addlegendentry{P1 elements (2.0)}
\addplot+[draw=none, color = {rgb,1:red,0.00000000;green,0.50196078;blue,0.00000000},
draw opacity = 1.0,
line width = 0,
solid,mark = triangle*,
mark size = 3.0,
mark options = {
    color = {rgb,1:red,0.00000000;green,0.00000000;blue,0.00000000}, draw opacity = 1.0,
    fill = {rgb,1:red,0.00000000;green,0.50196078;blue,0.00000000}, fill opacity = 1.0,
    line width = 1,
    rotate = 0,
    solid
},forget plot] coordinates {
(0.2776801836348979, 0.007899632948407665)
};
\addplot+[draw=none, color = {rgb,1:red,0.00000000;green,0.50196078;blue,0.00000000},
draw opacity = 1.0,
line width = 0,
solid,mark = triangle*,
mark size = 3.0,
mark options = {
    color = {rgb,1:red,0.00000000;green,0.00000000;blue,0.00000000}, draw opacity = 1.0,
    fill = {rgb,1:red,0.00000000;green,0.50196078;blue,0.00000000}, fill opacity = 1.0,
    line width = 1,
    rotate = 0,
    solid
},forget plot] coordinates {
(0.13884009181744894, 0.0020240250015253333)
};
\addplot+[draw=none, color = {rgb,1:red,0.00000000;green,0.50196078;blue,0.00000000},
draw opacity = 1.0,
line width = 0,
solid,mark = triangle*,
mark size = 3.0,
mark options = {
    color = {rgb,1:red,0.00000000;green,0.00000000;blue,0.00000000}, draw opacity = 1.0,
    fill = {rgb,1:red,0.00000000;green,0.50196078;blue,0.00000000}, fill opacity = 1.0,
    line width = 1,
    rotate = 0,
    solid
},forget plot] coordinates {
(0.06942004590872447, 0.0005094639292943957)
};
\addplot+[draw=none, color = {rgb,1:red,0.00000000;green,0.50196078;blue,0.00000000},
draw opacity = 1.0,
line width = 0,
solid,mark = triangle*,
mark size = 3.0,
mark options = {
    color = {rgb,1:red,0.00000000;green,0.00000000;blue,0.00000000}, draw opacity = 1.0,
    fill = {rgb,1:red,0.00000000;green,0.50196078;blue,0.00000000}, fill opacity = 1.0,
    line width = 1,
    rotate = 0,
    solid
},forget plot] coordinates {
(0.034710022954362235, 0.00012758990164437697)
};
\addplot+[draw=none, color = {rgb,1:red,0.00000000;green,0.50196078;blue,0.00000000},
draw opacity = 1.0,
line width = 0,
solid,mark = triangle*,
mark size = 3.0,
mark options = {
    color = {rgb,1:red,0.00000000;green,0.00000000;blue,0.00000000}, draw opacity = 1.0,
    fill = {rgb,1:red,0.00000000;green,0.50196078;blue,0.00000000}, fill opacity = 1.0,
    line width = 1,
    rotate = 0,
    solid
},forget plot] coordinates {
(0.017355011477181118, 3.191154208467998e-5)
};
\addplot+ [color = {rgb,1:red,0.00000000;green,0.50196078;blue,0.00000000},
draw opacity = 1.0,
line width = 1,
solid,mark = none,
mark size = 2.0,
mark options = {
    color = {rgb,1:red,0.00000000;green,0.00000000;blue,0.00000000}, draw opacity = 1.0,
    fill = {rgb,1:red,0.00000000;green,0.50196078;blue,0.00000000}, fill opacity = 1.0,
    line width = 1,
    rotate = 0,
    solid
},forget plot]coordinates {
(0.5553603672697958, 0.03267741909471229)
(0.2776801836348979, 0.00816935477367807)
(0.13884009181744894, 0.002042338693419519)
(0.06942004590872447, 0.0005105846733548797)
(0.034710022954362235, 0.00012764616833871987)
(0.017355011477181118, 3.1911542084679994e-5)
};
\addplot+ [color = {rgb,1:red,0.00000000;green,0.50196078;blue,0.00000000},
draw opacity = 1.0,
line width = 1,
solid,mark = none,
mark size = 2.0,
mark options = {
    color = {rgb,1:red,0.00000000;green,0.00000000;blue,0.00000000}, draw opacity = 1.0,
    fill = {rgb,1:red,0.00000000;green,0.50196078;blue,0.00000000}, fill opacity = 1.0,
    line width = 1,
    rotate = 0,
    solid
},forget plot]coordinates {
(0.5553603672697958, 0.0015440808887540907)
(0.2776801836348979, 0.00019301011109426126)
(0.13884009181744894, 2.4126263886782674e-5)
(0.06942004590872447, 3.015782985847833e-6)
(0.034710022954362235, 3.7697287323097896e-7)
(0.017355011477181118, 4.712160915387245e-8)
};
\end{axis}

\end{tikzpicture}
\end{subfigure}
\caption{Standard map: Errors in the first nontrivial eigenvalue (left) and 2-dimensional eigenspace corresponding to the two smallest nontrivial eigenvalues (right) of the dynamic Laplacian with the CG approach. The slopes of the corresponding lines are given in brackets in the legends.}
\label{fig:std_map_CG}
\end{figure}
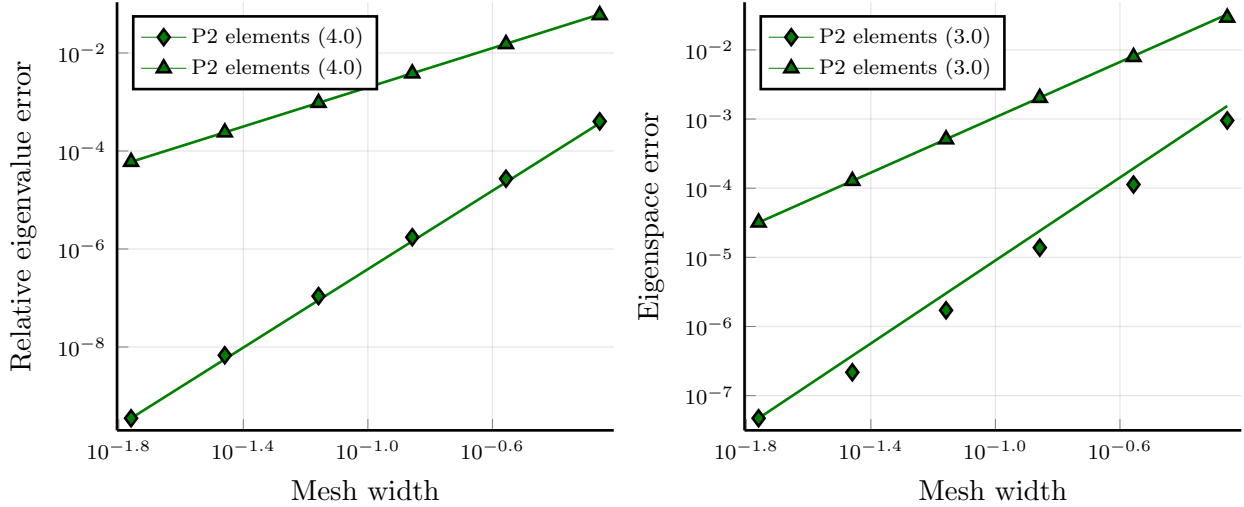

As a simple first example, we consider two iterations of the standard map
\begin{align*}
T(x,y) = (x + y + a\sin(x), y + a\sin(x))
\end{align*}
on the $2$-torus $\mathbb S^1 \times \mathbb S^1$ with parameter  $a = 0.971635$. This is the first example considered in \cite{rbf} and is a weakly nonlinear map.  Figure \ref{fig:std_map_CG} shows the eigenvalue and eigenspace errors in dependence of the mesh width for triangular $P^1$ and $P^2$ Lagrange elements.

The reference solution was computed with $P^2$ elements on a regular $1025\times1025$ grid with quadrature order $5$.
In this case, the experimentally observed rates almost perfectly agree with the predictions
\eqref{eq:CG_lambda_error} and \eqref{eq:CG_v_error_2} even for quadrature order $2$ in the computation of the entries of stiffness- and mass-matrix.
Improvements in errors ranged between one and two orders of magnitude when moving from $P^1$ to $P^2$ elements.
To be able to directly compare these results with the corresponding plots for TO methods in Figure \ref{fig:std_map_TO}, we used a two-dimensional eigenspace in Figure \ref{fig:std_map_CG}.

\subsubsection{Cylinder Flow}

As a second example, we consider the cylinder flow map used in \cite{rbf} based on \cite{froyland2010coherent}. This is a time-dependent flow on the cylinder\footnote{More precisely, we looked at the flow on $\mathbb S^1 \times [10^{-2}, \pi - 10^{-2}]$ as our numerical approximation of the flow-map was not well-behaved at the boundary }
$\mathbb S^1 \times [0,\pi]$ defined by the non-autonomous ordinary differential equation
\begin{align*}
\dot x(t) &= c - A(t)\sin(x - \nu t)\cos(y) + \varepsilon G(g(x,y,t))\sin(t/2) \\
\dot y(t) &= A(t)\cos(x - \nu t)\sin(y)
\end{align*}
where $A(t) = 1 + 0.125 \sin(2\sqrt 5 t)$, $G(\psi) = 1/(\psi^2 + 1)^2$ and $g(x,y,t) = \sin(x-\nu t)\sin(y) + y/2 - \pi/4$. Here the parameters $c=0.5, \nu=0.5,\varepsilon=0.25$ were used, the time-interval was $[0,40]$.

The reference solution was computed on a regular mesh on $1025\times 1025$ nodes with triangular $P^2$-Lagrange elements and quadrature order $8$. We used quadrature order $5$ for the numerical experiments (since we did not observe the same rates for smaller orders).

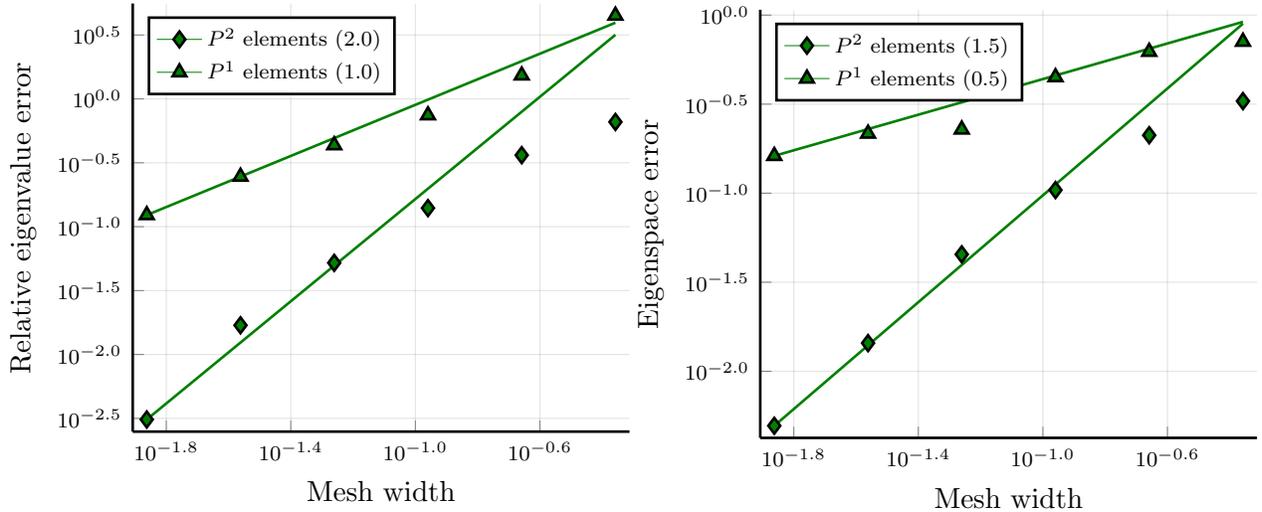
\begin{figure}[H]
\captionsetup[subfigure]{justification=centering}
\begin{subfigure}{0.5\linewidth}
\begin{tikzpicture}[]
\begin{axis}[scale only axis,width=0.8\linewidth, legend pos = {north west}, ylabel = {Relative eigenvalue error}, xmin = {0.012365463918175449}, xmax = {0.487157498280955}, ymax = {5.579215932300088}, ymode = {log}, xlabel = {Mesh width}, unbounded coords=jump,scaled x ticks = false,xlabel style = {font = {\fontsize{11 pt}{14.3 pt}\selectfont}, color = {rgb,1:red,0.00000000;green,0.00000000;blue,0.00000000}, draw opacity = 1.0, rotate = 0.0},log basis x=10,xmajorgrids = true,xtick = {0.015848931924611134,0.039810717055349734,0.1,0.251188643150958},xticklabels = {$10^{-1.8}$,$10^{-1.4}$,$10^{-1.0}$,$10^{-0.6}$},xtick align = inside,xticklabel style = {font = {\fontsize{8 pt}{10.4 pt}\selectfont}, color = {rgb,1:red,0.00000000;green,0.00000000;blue,0.00000000}, draw opacity = 1.0, rotate = 0.0},x grid style = {color = {rgb,1:red,0.00000000;green,0.00000000;blue,0.00000000},
draw opacity = 0.1,
line width = 0.5,
solid},axis lines* = left,x axis line style = {color = {rgb,1:red,0.00000000;green,0.00000000;blue,0.00000000},
draw opacity = 1.0,
line width = 1,
solid},scaled y ticks = false,ylabel style = {font = {\fontsize{11 pt}{14.3 pt}\selectfont}, color = {rgb,1:red,0.00000000;green,0.00000000;blue,0.00000000}, draw opacity = 1.0, rotate = 0.0},log basis y=10,ymajorgrids = true,ytick = {0.0031622776601683794,0.01,0.03162277660168379,0.1,0.31622776601683794,1.0,3.162277660168379},yticklabels = {$10^{-2.5}$,$10^{-2.0}$,$10^{-1.5}$,$10^{-1.0}$,$10^{-0.5}$,$10^{0.0}$,$10^{0.5}$},ytick align = inside,yticklabel style = {font = {\fontsize{8 pt}{10.4 pt}\selectfont}, color = {rgb,1:red,0.00000000;green,0.00000000;blue,0.00000000}, draw opacity = 1.0, rotate = 0.0},y grid style = {color = {rgb,1:red,0.00000000;green,0.00000000;blue,0.00000000},
draw opacity = 0.1,
line width = 0.5,
solid},axis lines* = left,y axis line style = {color = {rgb,1:red,0.00000000;green,0.00000000;blue,0.00000000},
draw opacity = 1.0,
line width = 1,
solid},    xshift = 0.0mm,
    yshift = 0.0mm,
    axis background/.style={fill={rgb,1:red,1.00000000;green,1.00000000;blue,1.00000000}}
,legend style = {color = {rgb,1:red,0.00000000;green,0.00000000;blue,0.00000000},
draw opacity = 1.0,
line width = 1,
solid,fill = {rgb,1:red,1.00000000;green,1.00000000;blue,1.00000000},font = {\fontsize{8 pt}{10.4 pt}\selectfont}},colorbar style={title=}, xmode = {log}, ymin = {0.002491080221658315}]\addplot+[draw=none, color = {rgb,1:red,0.00000000;green,0.50196078;blue,0.00000000},
draw opacity = 1.0,
line width = 0,
solid,mark = diamond*,
mark size = 3.0,
mark options = {
    color = {rgb,1:red,0.00000000;green,0.00000000;blue,0.00000000}, draw opacity = 1.0,
    fill = {rgb,1:red,0.00000000;green,0.50196078;blue,0.00000000}, fill opacity = 1.0,
    line width = 1,
    rotate = 0,
    solid
}] coordinates {
(0.4390509206900454, 0.6599608903257889)
};
\addlegendentry{$P^2$ elements (2.0)}
\addplot+[draw=none, color = {rgb,1:red,0.00000000;green,0.50196078;blue,0.00000000},
draw opacity = 1.0,
line width = 0,
solid,mark = diamond*,
mark size = 3.0,
mark options = {
    color = {rgb,1:red,0.00000000;green,0.00000000;blue,0.00000000}, draw opacity = 1.0,
    fill = {rgb,1:red,0.00000000;green,0.50196078;blue,0.00000000}, fill opacity = 1.0,
    line width = 1,
    rotate = 0,
    solid
},forget plot] coordinates {
(0.2195254603450227, 0.36290598466646246)
};
\addplot+[draw=none, color = {rgb,1:red,0.00000000;green,0.50196078;blue,0.00000000},
draw opacity = 1.0,
line width = 0,
solid,mark = diamond*,
mark size = 3.0,
mark options = {
    color = {rgb,1:red,0.00000000;green,0.00000000;blue,0.00000000}, draw opacity = 1.0,
    fill = {rgb,1:red,0.00000000;green,0.50196078;blue,0.00000000}, fill opacity = 1.0,
    line width = 1,
    rotate = 0,
    solid
},forget plot] coordinates {
(0.10976273017251136, 0.13982535143316713)
};
\addplot+[draw=none, color = {rgb,1:red,0.00000000;green,0.50196078;blue,0.00000000},
draw opacity = 1.0,
line width = 0,
solid,mark = diamond*,
mark size = 3.0,
mark options = {
    color = {rgb,1:red,0.00000000;green,0.00000000;blue,0.00000000}, draw opacity = 1.0,
    fill = {rgb,1:red,0.00000000;green,0.50196078;blue,0.00000000}, fill opacity = 1.0,
    line width = 1,
    rotate = 0,
    solid
},forget plot] coordinates {
(0.05488136508625568, 0.05210905117124607)
};
\addplot+[draw=none, color = {rgb,1:red,0.00000000;green,0.50196078;blue,0.00000000},
draw opacity = 1.0,
line width = 0,
solid,mark = diamond*,
mark size = 3.0,
mark options = {
    color = {rgb,1:red,0.00000000;green,0.00000000;blue,0.00000000}, draw opacity = 1.0,
    fill = {rgb,1:red,0.00000000;green,0.50196078;blue,0.00000000}, fill opacity = 1.0,
    line width = 1,
    rotate = 0,
    solid
},forget plot] coordinates {
(0.02744068254312784, 0.016923244509504887)
};
\addplot+[draw=none, color = {rgb,1:red,0.00000000;green,0.50196078;blue,0.00000000},
draw opacity = 1.0,
line width = 0,
solid,mark = diamond*,
mark size = 3.0,
mark options = {
    color = {rgb,1:red,0.00000000;green,0.00000000;blue,0.00000000}, draw opacity = 1.0,
    fill = {rgb,1:red,0.00000000;green,0.50196078;blue,0.00000000}, fill opacity = 1.0,
    line width = 1,
    rotate = 0,
    solid
},forget plot] coordinates {
(0.01372034127156392, 0.0030988766037429527)
};
\addplot+[draw=none, color = {rgb,1:red,0.00000000;green,0.50196078;blue,0.00000000},
draw opacity = 1.0,
line width = 0,
solid,mark = triangle*,
mark size = 3.0,
mark options = {
    color = {rgb,1:red,0.00000000;green,0.00000000;blue,0.00000000}, draw opacity = 1.0,
    fill = {rgb,1:red,0.00000000;green,0.50196078;blue,0.00000000}, fill opacity = 1.0,
    line width = 1,
    rotate = 0,
    solid
}] coordinates {
(0.4390509206900454, 4.484939621192657)
};
\addlegendentry{$P^1$ elements (1.0)}
\addplot+[draw=none, color = {rgb,1:red,0.00000000;green,0.50196078;blue,0.00000000},
draw opacity = 1.0,
line width = 0,
solid,mark = triangle*,
mark size = 3.0,
mark options = {
    color = {rgb,1:red,0.00000000;green,0.00000000;blue,0.00000000}, draw opacity = 1.0,
    fill = {rgb,1:red,0.00000000;green,0.50196078;blue,0.00000000}, fill opacity = 1.0,
    line width = 1,
    rotate = 0,
    solid
},forget plot] coordinates {
(0.2195254603450227, 1.5266835690246403)
};
\addplot+[draw=none, color = {rgb,1:red,0.00000000;green,0.50196078;blue,0.00000000},
draw opacity = 1.0,
line width = 0,
solid,mark = triangle*,
mark size = 3.0,
mark options = {
    color = {rgb,1:red,0.00000000;green,0.00000000;blue,0.00000000}, draw opacity = 1.0,
    fill = {rgb,1:red,0.00000000;green,0.50196078;blue,0.00000000}, fill opacity = 1.0,
    line width = 1,
    rotate = 0,
    solid
},forget plot] coordinates {
(0.10976273017251136, 0.7448341709586542)
};
\addplot+[draw=none, color = {rgb,1:red,0.00000000;green,0.50196078;blue,0.00000000},
draw opacity = 1.0,
line width = 0,
solid,mark = triangle*,
mark size = 3.0,
mark options = {
    color = {rgb,1:red,0.00000000;green,0.00000000;blue,0.00000000}, draw opacity = 1.0,
    fill = {rgb,1:red,0.00000000;green,0.50196078;blue,0.00000000}, fill opacity = 1.0,
    line width = 1,
    rotate = 0,
    solid
},forget plot] coordinates {
(0.05488136508625568, 0.434096576852428)
};
\addplot+[draw=none, color = {rgb,1:red,0.00000000;green,0.50196078;blue,0.00000000},
draw opacity = 1.0,
line width = 0,
solid,mark = triangle*,
mark size = 3.0,
mark options = {
    color = {rgb,1:red,0.00000000;green,0.00000000;blue,0.00000000}, draw opacity = 1.0,
    fill = {rgb,1:red,0.00000000;green,0.50196078;blue,0.00000000}, fill opacity = 1.0,
    line width = 1,
    rotate = 0,
    solid
},forget plot] coordinates {
(0.02744068254312784, 0.24627988705429862)
};
\addplot+[draw=none, color = {rgb,1:red,0.00000000;green,0.50196078;blue,0.00000000},
draw opacity = 1.0,
line width = 0,
solid,mark = triangle*,
mark size = 3.0,
mark options = {
    color = {rgb,1:red,0.00000000;green,0.00000000;blue,0.00000000}, draw opacity = 1.0,
    fill = {rgb,1:red,0.00000000;green,0.50196078;blue,0.00000000}, fill opacity = 1.0,
    line width = 1,
    rotate = 0,
    solid
},forget plot] coordinates {
(0.01372034127156392, 0.12329416803269028)
};
\addplot+ [color = {rgb,1:red,0.00000000;green,0.50196078;blue,0.00000000},
draw opacity = 1.0,
line width = 1,
solid,mark = none,
mark size = 2.0,
mark options = {
    color = {rgb,1:red,0.00000000;green,0.00000000;blue,0.00000000}, draw opacity = 1.0,
    fill = {rgb,1:red,0.00000000;green,0.50196078;blue,0.00000000}, fill opacity = 1.0,
    line width = 1,
    rotate = 0,
    solid
},forget plot]coordinates {
(0.4390509206900454, 3.9454133770460897)
(0.2195254603450227, 1.9727066885230453)
(0.10976273017251136, 0.9863533442615224)
(0.05488136508625568, 0.4931766721307613)
(0.02744068254312784, 0.2465883360653806)
(0.01372034127156392, 0.12329416803269029)
};
\addplot+ [color = {rgb,1:red,0.00000000;green,0.50196078;blue,0.00000000},
draw opacity = 1.0,
line width = 1,
solid,mark = none,
mark size = 2.0,
mark options = {
    color = {rgb,1:red,0.00000000;green,0.00000000;blue,0.00000000}, draw opacity = 1.0,
    fill = {rgb,1:red,0.00000000;green,0.50196078;blue,0.00000000}, fill opacity = 1.0,
    line width = 1,
    rotate = 0,
    solid
},forget plot]coordinates {
(0.4390509206900454, 3.1732496422327836)
(0.2195254603450227, 0.7933124105581961)
(0.10976273017251136, 0.19832810263954898)
(0.05488136508625568, 0.04958202565988726)
(0.02744068254312784, 0.012395506414971811)
(0.01372034127156392, 0.003098876603742952)
};
\end{axis}

\end{tikzpicture}
\end{subfigure}
\begin{subfigure}{0.5\linewidth}
\begin{tikzpicture}[]
\begin{axis}[scale only axis,width=0.8\linewidth, legend pos = {north west}, ylabel = {Eigenspace  error}, xmin = {0.012365463918175449}, xmax = {0.487157498280955}, ymax = {1.0714608895408926}, ymode = {log}, xlabel = {Mesh width}, unbounded coords=jump,scaled x ticks = false,xlabel style = {font = {\fontsize{11 pt}{14.3 pt}\selectfont}, color = {rgb,1:red,0.00000000;green,0.00000000;blue,0.00000000}, draw opacity = 1.0, rotate = 0.0},log basis x=10,xmajorgrids = true,xtick = {0.015848931924611134,0.039810717055349734,0.1,0.251188643150958},xticklabels = {$10^{-1.8}$,$10^{-1.4}$,$10^{-1.0}$,$10^{-0.6}$},xtick align = inside,xticklabel style = {font = {\fontsize{8 pt}{10.4 pt}\selectfont}, color = {rgb,1:red,0.00000000;green,0.00000000;blue,0.00000000}, draw opacity = 1.0, rotate = 0.0},x grid style = {color = {rgb,1:red,0.00000000;green,0.00000000;blue,0.00000000},
draw opacity = 0.1,
line width = 0.5,
solid},axis lines* = left,x axis line style = {color = {rgb,1:red,0.00000000;green,0.00000000;blue,0.00000000},
draw opacity = 1.0,
line width = 1,
solid},scaled y ticks = false,ylabel style = {font = {\fontsize{11 pt}{14.3 pt}\selectfont}, color = {rgb,1:red,0.00000000;green,0.00000000;blue,0.00000000}, draw opacity = 1.0, rotate = 0.0},log basis y=10,ymajorgrids = true,ytick = {0.01,0.03162277660168379,0.1,0.31622776601683794,1.0},yticklabels = {$10^{-2.0}$,$10^{-1.5}$,$10^{-1.0}$,$10^{-0.5}$,$10^{0.0}$},ytick align = inside,yticklabel style = {font = {\fontsize{8 pt}{10.4 pt}\selectfont}, color = {rgb,1:red,0.00000000;green,0.00000000;blue,0.00000000}, draw opacity = 1.0, rotate = 0.0},y grid style = {color = {rgb,1:red,0.00000000;green,0.00000000;blue,0.00000000},
draw opacity = 0.1,
line width = 0.5,
solid},axis lines* = left,y axis line style = {color = {rgb,1:red,0.00000000;green,0.00000000;blue,0.00000000},
draw opacity = 1.0,
line width = 1,
solid},    xshift = 0.0mm,
    yshift = 0.0mm,
    axis background/.style={fill={rgb,1:red,1.00000000;green,1.00000000;blue,1.00000000}}
,legend style = {color = {rgb,1:red,0.00000000;green,0.00000000;blue,0.00000000},
draw opacity = 1.0,
line width = 1,
solid,fill = {rgb,1:red,1.00000000;green,1.00000000;blue,1.00000000},font = {\fontsize{8 pt}{10.4 pt}\selectfont}},colorbar style={title=}, xmode = {log}, ymin = {0.004226880900655178}]\addplot+[draw=none, color = {rgb,1:red,0.00000000;green,0.50196078;blue,0.00000000},
draw opacity = 1.0,
line width = 0,
solid,mark = diamond*,
mark size = 3.0,
mark options = {
    color = {rgb,1:red,0.00000000;green,0.00000000;blue,0.00000000}, draw opacity = 1.0,
    fill = {rgb,1:red,0.00000000;green,0.50196078;blue,0.00000000}, fill opacity = 1.0,
    line width = 1,
    rotate = 0,
    solid
}] coordinates {
(0.4390509206900454, 0.328884826119691)
};
\addlegendentry{$P^2$ elements (1.5)}
\addplot+[draw=none, color = {rgb,1:red,0.00000000;green,0.50196078;blue,0.00000000},
draw opacity = 1.0,
line width = 0,
solid,mark = diamond*,
mark size = 3.0,
mark options = {
    color = {rgb,1:red,0.00000000;green,0.00000000;blue,0.00000000}, draw opacity = 1.0,
    fill = {rgb,1:red,0.00000000;green,0.50196078;blue,0.00000000}, fill opacity = 1.0,
    line width = 1,
    rotate = 0,
    solid
},forget plot] coordinates {
(0.2195254603450227, 0.21157712846425608)
};
\addplot+[draw=none, color = {rgb,1:red,0.00000000;green,0.50196078;blue,0.00000000},
draw opacity = 1.0,
line width = 0,
solid,mark = diamond*,
mark size = 3.0,
mark options = {
    color = {rgb,1:red,0.00000000;green,0.00000000;blue,0.00000000}, draw opacity = 1.0,
    fill = {rgb,1:red,0.00000000;green,0.50196078;blue,0.00000000}, fill opacity = 1.0,
    line width = 1,
    rotate = 0,
    solid
},forget plot] coordinates {
(0.10976273017251136, 0.10413829272830695)
};
\addplot+[draw=none, color = {rgb,1:red,0.00000000;green,0.50196078;blue,0.00000000},
draw opacity = 1.0,
line width = 0,
solid,mark = diamond*,
mark size = 3.0,
mark options = {
    color = {rgb,1:red,0.00000000;green,0.00000000;blue,0.00000000}, draw opacity = 1.0,
    fill = {rgb,1:red,0.00000000;green,0.50196078;blue,0.00000000}, fill opacity = 1.0,
    line width = 1,
    rotate = 0,
    solid
},forget plot] coordinates {
(0.05488136508625568, 0.04534711988338571)
};
\addplot+[draw=none, color = {rgb,1:red,0.00000000;green,0.50196078;blue,0.00000000},
draw opacity = 1.0,
line width = 0,
solid,mark = diamond*,
mark size = 3.0,
mark options = {
    color = {rgb,1:red,0.00000000;green,0.00000000;blue,0.00000000}, draw opacity = 1.0,
    fill = {rgb,1:red,0.00000000;green,0.50196078;blue,0.00000000}, fill opacity = 1.0,
    line width = 1,
    rotate = 0,
    solid
},forget plot] coordinates {
(0.02744068254312784, 0.014393424967977542)
};
\addplot+[draw=none, color = {rgb,1:red,0.00000000;green,0.50196078;blue,0.00000000},
draw opacity = 1.0,
line width = 0,
solid,mark = diamond*,
mark size = 3.0,
mark options = {
    color = {rgb,1:red,0.00000000;green,0.00000000;blue,0.00000000}, draw opacity = 1.0,
    fill = {rgb,1:red,0.00000000;green,0.50196078;blue,0.00000000}, fill opacity = 1.0,
    line width = 1,
    rotate = 0,
    solid
},forget plot] coordinates {
(0.01372034127156392, 0.004943750109196259)
};
\addplot+[draw=none, color = {rgb,1:red,0.00000000;green,0.50196078;blue,0.00000000},
draw opacity = 1.0,
line width = 0,
solid,mark = triangle*,
mark size = 3.0,
mark options = {
    color = {rgb,1:red,0.00000000;green,0.00000000;blue,0.00000000}, draw opacity = 1.0,
    fill = {rgb,1:red,0.00000000;green,0.50196078;blue,0.00000000}, fill opacity = 1.0,
    line width = 1,
    rotate = 0,
    solid
}] coordinates {
(0.4390509206900454, 0.7088903825840833)
};
\addlegendentry{$P^1$ elements (0.5)}
\addplot+[draw=none, color = {rgb,1:red,0.00000000;green,0.50196078;blue,0.00000000},
draw opacity = 1.0,
line width = 0,
solid,mark = triangle*,
mark size = 3.0,
mark options = {
    color = {rgb,1:red,0.00000000;green,0.00000000;blue,0.00000000}, draw opacity = 1.0,
    fill = {rgb,1:red,0.00000000;green,0.50196078;blue,0.00000000}, fill opacity = 1.0,
    line width = 1,
    rotate = 0,
    solid
},forget plot] coordinates {
(0.2195254603450227, 0.6223750093023916)
};
\addplot+[draw=none, color = {rgb,1:red,0.00000000;green,0.50196078;blue,0.00000000},
draw opacity = 1.0,
line width = 0,
solid,mark = triangle*,
mark size = 3.0,
mark options = {
    color = {rgb,1:red,0.00000000;green,0.00000000;blue,0.00000000}, draw opacity = 1.0,
    fill = {rgb,1:red,0.00000000;green,0.50196078;blue,0.00000000}, fill opacity = 1.0,
    line width = 1,
    rotate = 0,
    solid
},forget plot] coordinates {
(0.10976273017251136, 0.44819991513993745)
};
\addplot+[draw=none, color = {rgb,1:red,0.00000000;green,0.50196078;blue,0.00000000},
draw opacity = 1.0,
line width = 0,
solid,mark = triangle*,
mark size = 3.0,
mark options = {
    color = {rgb,1:red,0.00000000;green,0.00000000;blue,0.00000000}, draw opacity = 1.0,
    fill = {rgb,1:red,0.00000000;green,0.50196078;blue,0.00000000}, fill opacity = 1.0,
    line width = 1,
    rotate = 0,
    solid
},forget plot] coordinates {
(0.05488136508625568, 0.22783078515234292)
};
\addplot+[draw=none, color = {rgb,1:red,0.00000000;green,0.50196078;blue,0.00000000},
draw opacity = 1.0,
line width = 0,
solid,mark = triangle*,
mark size = 3.0,
mark options = {
    color = {rgb,1:red,0.00000000;green,0.00000000;blue,0.00000000}, draw opacity = 1.0,
    fill = {rgb,1:red,0.00000000;green,0.50196078;blue,0.00000000}, fill opacity = 1.0,
    line width = 1,
    rotate = 0,
    solid
},forget plot] coordinates {
(0.02744068254312784, 0.21601735551649853)
};
\addplot+[draw=none, color = {rgb,1:red,0.00000000;green,0.50196078;blue,0.00000000},
draw opacity = 1.0,
line width = 0,
solid,mark = triangle*,
mark size = 3.0,
mark options = {
    color = {rgb,1:red,0.00000000;green,0.00000000;blue,0.00000000}, draw opacity = 1.0,
    fill = {rgb,1:red,0.00000000;green,0.50196078;blue,0.00000000}, fill opacity = 1.0,
    line width = 1,
    rotate = 0,
    solid
},forget plot] coordinates {
(0.01372034127156392, 0.16194398970624382)
};
\addplot+ [color = {rgb,1:red,0.00000000;green,0.50196078;blue,0.00000000},
draw opacity = 1.0,
line width = 1,
solid,mark = none,
mark size = 2.0,
mark options = {
    color = {rgb,1:red,0.00000000;green,0.00000000;blue,0.00000000}, draw opacity = 1.0,
    fill = {rgb,1:red,0.00000000;green,0.50196078;blue,0.00000000}, fill opacity = 1.0,
    line width = 1,
    rotate = 0,
    solid
},forget plot]coordinates {
(0.4390509206900454, 0.9160935463495156)
(0.2195254603450227, 0.6477759588249753)
(0.10976273017251136, 0.4580467731747578)
(0.05488136508625568, 0.32388797941248765)
(0.02744068254312784, 0.2290233865873789)
(0.01372034127156392, 0.1619439897062438)
};
\addplot+ [color = {rgb,1:red,0.00000000;green,0.50196078;blue,0.00000000},
draw opacity = 1.0,
line width = 1,
solid,mark = none,
mark size = 2.0,
mark options = {
    color = {rgb,1:red,0.00000000;green,0.00000000;blue,0.00000000}, draw opacity = 1.0,
    fill = {rgb,1:red,0.00000000;green,0.50196078;blue,0.00000000}, fill opacity = 1.0,
    line width = 1,
    rotate = 0,
    solid
},forget plot]coordinates {
(0.4390509206900454, 0.8949143620363285)
(0.2195254603450227, 0.3164000069885605)
(0.10976273017251136, 0.11186429525454107)
(0.05488136508625568, 0.03955000087357007)
(0.02744068254312784, 0.013983036906817636)
(0.01372034127156392, 0.004943750109196256)
};
\end{axis}

\end{tikzpicture}
\end{subfigure}
 \caption{Cylinder flow: Errors in the first nontrivial eigenvalue (left) and corresponding 1-dimensional eigenspace (right) of the dynamic Laplacian discretized with CG approach. The slopes of the corresponding lines are given in brackets in the legends.}
 \label{fig:cyl_CG}
\end{figure}

The orders of convergence observed for the eigenvalue and eigenvector errors is surprisingly low;  the experimental values are almost exactly half of the orders predicted in \eqref{eq:CG_lambda_error}, \eqref{eq:CG_v_error_1} or \eqref{eq:CG_v_error_2}.
The observed order of the eigenvector errors for $P^1$ elements is particularly poor.
We also note that the slopes shown in the figure do not remain consistent when varying the quadrature order.

For this flow (and other highly stretching and highly nonlinear systems), the mean diffusion tensor $A = [DT_t]^{-1}([DT_t]^{-1})^T$ used in the weak formulation of the dynamic Laplacian (\ref{eq:weak_form_CG}) has extremely high variation locally (cf.~Figure \ref{fig:cylinder_dbs}).  We suspect that even the reference grid of $1025\times 1025$ nodes is not sufficiently resolved to be in a regime where the convergence orders predicted by the theory can be observed.

\begin{figure}[H]
\centering
\input{Figures/processed/cylinder_dbs.tex}
\caption{Cylinder flow: log of the trace of the mean diffusion tensor, $\log_{10}\mbox{tr}(A(x))$.}
\label{fig:cylinder_dbs}
\end{figure}

\subsection{Bickley jet}
\label{eq:exp_CG_Bickley}

\begin{figure}[H]
\centering
\begin{subfigure}{\linewidth}
\centering
\input{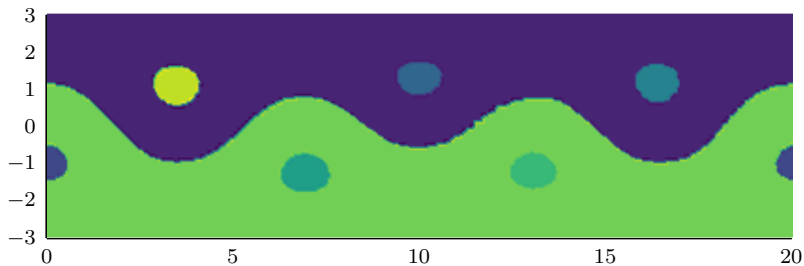}
\subcaption{$P^1$ elements on $101\times 31$ mesh (18000 quadrature points, 6000 triangles)}
\end{subfigure}
\begin{subfigure}{\linewidth}
\centering
\input{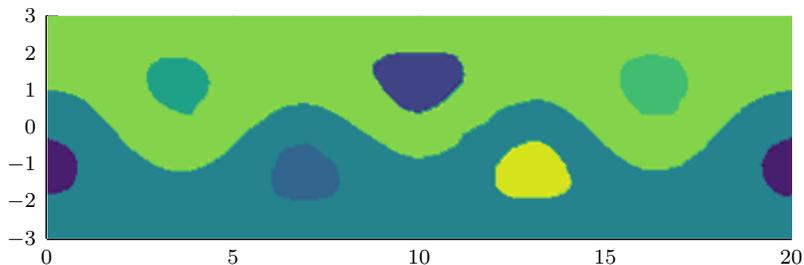}
\subcaption{$P^2$ elements on $21\times 7$ mesh (720 quadrature points, 240 triangles)}
\end{subfigure}
\caption{Bickley jet: Comparison of coherent sets obtained by the CG approach with $P^1$ vs.\ $P^2$ elements.}
\label{fig:bickley_cg}
\end{figure}

The Bickley jet flow is a well-known test case first introduced in \cite{rypina2007lagrangian}. The flow is defined by a stream-function
\begin{align}
\psi(x,y,t) &= -U_0L \tanh(y/L_0) + \sum_{i=1}^3 A_i U_0 L \mbox{sech}^2(y/L) \cos(k_i(x - c_i t))
\end{align}
with constants $U_0,L_0, A_i, k_i,c_i$ on a cylindrical domain (see \cite{rypina2007lagrangian}). We considered the flow for a timespan of $40$ days.
$8$-partitions were computed by k-means clustering on $200\times60$ values of the leading eigenvectors of the FEM approximation.
The mesh widths in Figure \ref{fig:bickley_cg} below are the lowest (at the aspect ratio $10:3$) for which the topology of the clustering result does not change. It is possible to obtain the same result (topologically)  with $P^2$ elements at significantly
reduced cost (in terms of the number of quadrature point used) to $P^1$ elements. In many cases, the computation requiring the 
most runtime is the computation of diffusion tensors by time-integration at the quadrature points. 
As shown in Figure \ref{fig:bickley_cg}, in this case this is reduced by a factor of 25 by using $P^2$ elements.

\par

\section{Convergence for the TO approach}

The CG approach has the disadvantage of requiring the numerical approximation of the derivative $DT_t$ and a subsequent quadrature which might be challenging if $DT_t$ has high variation locally.

An alternative approach \cite{froylandjunge18} is to rewrite the weak form (\ref{eq:weak_form_CG}) of the dynamic Laplacian operator using the transfer operator $T_t$. 
In the sequel we discuss the case $t \in \mathcal{I} := \{0,1\}$, although the results hold for general finite $\mathcal{I}$.
Letting $a^t$ be the weak-form of the Laplacian on $\Omega_t$, we have
\begin{align*}
a(u,v) = \frac{1}{2}\left(a^0(u,v) + a^1(Tu, Tv)\right)
\end{align*}
where we write $T := T_{1}$, and for brevity overload this notation so that $T$ acts on elements of $\Omega_0$ in the usual way \emph{and} acts on functions in $L^2(\Omega_0)$ as $Tu:=T_{1,*}u$.
In the transfer operator (TO) approach, we approximate $T$ on a suitable subspace of $H^1$ by
an operator of the form $I_h T$, where $I_h$ is some projection operator.
Let $\cT^1_h$ be a quasi-uniform family of triangulations on $\Omega_1$,
$S^1_h$ the approximation space with Lagrange $P^{k'}$ finite elements and $\hat S_h^1 := S_h^1 \cap \hat H^1(\Omega_1)$.
Note that in general, the collection $\cT_h^0$ (introduced at the beginning of Section \ref{sec:CGconv}) and $\cT^1_h$ are unrelated. Additionally, we assume that the meshes used are nested
in the sense that $S_{h'} \subset S_{h}$ for $h < h'$.

We consider the following options for $I_h$ in the sequel:
\begin{enumerate}
\item the $L^2$-orthogonal projection onto $S_h^1$ (\emph{$L^2$-Galerkin}),
\item the $H^1$-orthogonal projection onto $S_h^1$ (\emph{$H^1$-Galerkin}),
\item the (Lagrange $P^k$) nodal interpolation on $S_h^1$ (\emph{collocation}).
\end{enumerate}
Case 3.\ corresponds to the ``TO approach'' from \cite{froylandjunge18}, with the difference between the adaptive and non-adaptive approaches being the choice of the mesh $\cT^1_h$.

We now define the ($h$-dependent) symmetric bilinear form
\begin{align}\label{eq:ah}
\tilde a_h(u,v) := \frac{1}{2}\left( a^0(u,v) + a^1(I_h Tu, I_h Tv) \right)
\end{align}

As done in \cite{froylandjunge18}, eigenpairs are numerically approximated by calculating generalized eigenpairs of $(\tilde D_h,M_h)$ given by
\begin{align*}
(\tilde D_h)_{ij} &:= \tilde{a}_h(\varphi^{i}_h, \varphi^{j}_h)\\
(M_h)_{ij} &= \Bra{\varphi^{i}_h,\varphi^{j}_h}_{L^2}
\end{align*}

As $a^0$ is elliptic on $\hat H^1(\Omega_0)$, it follows immediately from \eqref{eq:ah} that $\tilde a_h$ is elliptic uniformly in $h$ on $\hat S_h^0$, with lower bound of the ellipticity constant being half of that of $a^0$.

\paragraph{Properties of $I_h$.} In order to investigate the convergence of the numerically calculated eigenpairs to those of the dynamic Laplacian, we will need to consider the following conditions for $I_h T$ acting on $S_h^0$:
\begin{enumerate}[(i)]
\item exactness for constant functions: $I_h T\chi_{\Omega_0} = \chi_{\Omega_1} $,
\item $H^1$-convergence: $\lim\limits_{h \rightarrow 0} \norm{I_h T f - Tf}_{H^1} = 0$ for any $f \in \cup_h S_h^0$,
\item $L^2$-stability:  $\sup_{h<h'}\norm{I_h T| S_{h'}^0}_{L^2 \rightarrow L^2} \leq C$.
\end{enumerate}
We now consider whether these three conditions are satisfied by the three variants of $I_h$:
\begin{enumerate}[(i)]
    \item It is easy to see that condition (i) is satisfied for all three variants, by the fact that $T := T_1$ is volume-preserving and $I_h$ is exact on constant functions.
    \item Condition (ii) is satisfied by the $H^1$-Galerkin approach. Under the technical assumption that $\mathcal T^1_h$ is quasi-uniform \cite{ernguermond}, we also get (ii) for the $L^2$-Galerkin approach.
    \item Condition (iii) is trivially satisfied by the $L^2$-Galerkin approach.
\end{enumerate}
We do not consider condition (iii) for the $H^1$-Galerkin approach. It remains to show (ii) and (iii) for the case of $I_h$ being the nodal interpolation operator.  In the appendix, we give  proofs of the following two Lemmas for $P^1$-triangular elements under the assumption that both meshes are quasi-uniform.
For $h'>h>0$, let $\cT^0_{h'}$ and $\cT^1_h$ be families of quasi-uniform triangulations of $\Omega_0$ and $\Omega_1$ respectively. Let $I_h : C(\overline{\Omega_1}) \rightarrow S_h^1$ be the $P^1$-Lagrange nodal interpolation operator corresponding to $\cT^1_h$.

\begin{lemma}[$L^2$-stability of $P^1$-nodal interpolation]
\label{lemma:stability}
There exists a constant $C > 0$ so that $\norm{I_h Tv}_{L^2} \leq C \norm{v}_{L^2}$ for all $v \in S_{h'}^0$ and all $h<h'$.
\end{lemma}

\begin{lemma}[$H^1$-convergence of $P^1$-nodal interpolation]
\label{lemma:convergence}
    If $v$ is piecewise $C^\infty$ on $\Omega_1$ then
    \[
    \norm{I_h v - v}_{H^1} \rightarrow 0 \qquad \text{as } h\to 0.
    \]
\end{lemma}

\paragraph{Convergence for homogeneous Neumann Boundary.}

Let $r(v)=a(v,v)/\Bra{v,v}_{L^2}$ and $\tilde{r}_h(v) = \tilde a_h(v,v)/\Bra{v,v}_{L^2}$ be the Rayleigh quotients for $a$ and $\tilde a_h$ respectively. Let $\lambda_1$ and $\tilde \lambda_{h,1}$ be the corresponding smallest Ritz values on $\hat H^1$ and $\hat S_h^0$, respectively and $v$, $ \tilde v_{h}$ the corresponding Ritz vectors, normalized to have unit $L^2$-norm.

From the aforementioned properties of $I_h T$ and the fact that the Ritz vectors minimize the Rayleigh quotient we now
derive a convergence result for $\tilde v_h$ to $v$ in two parts:
\begin{theorem}\label{theorem:upperbound}
Assume that $I_h$ satisfies the $H^1$-convergence property (ii), then
\begin{align*}
    \limsup\limits_{h \rightarrow 0} \tilde \lambda_{h,1} \leq \lambda_1
\end{align*}
\end{theorem}
\begin{proof}
    Pick any sequence $h_n$ with $h_n \rightarrow 0$. Then for arbitrary $w \in \hat S_h^0$
    \begin{align*}
        \limsup\limits_{n \rightarrow \infty} \tilde \lambda_{h_n,1} &= \limsup\limits_{n \rightarrow \infty}\tilde r_{h_n}(\tilde v_{h_n}) \ \ \mbox{(definition of $\tilde v_{h_n}$)} \\
        &\leq \limsup \limits_{n \rightarrow \infty} \tilde r_{h_n}(w)\ \  \mbox{ (Ritz-values minimize the Rayleigh-quotient)}
    \end{align*}
    Now for arbitrary $\epsilon > 0$, pick\footnote{That we can do this follows from the classical theory of convergence of eigenvalues of the
    Galerkin approximation such as equation \eqref{eq:CG_lambda_error} in the case that the eigenvectors are sufficiently smooth at the boundary.}
    some $w \in \hat S_{h'}^0$ where $h' = h'(\epsilon)$ with $\norm{w}_{L^2} = 1$ so that $r(w) < \lambda_1 + \epsilon$.
    Then by the assumption that our meshes are nested, $w$ is in the domain of $\tilde a_{h_n}$ for sufficiently large $n$ and we have:
    \begin{align*}
        \limsup\limits_{n \rightarrow \infty} \tilde \lambda_{h_n,1} &\leq \limsup\limits_{n \rightarrow \infty} \tilde{r}_{h_n}(w)&  \\
        &= \limsup\limits_{n \rightarrow \infty} \tilde a_{h_n}(w,w) \ \ \ \ \mbox{(definition of Rayleigh Quotient)}\\
        &= \limsup\limits_{n \rightarrow \infty}( a(w,w) + \tilde{a}_{h_n}(w,w) - a(w,w)) & \\
        &\leq \limsup\limits_{n \rightarrow \infty}\left( \lambda_1 + \epsilon + a^1(I_h Tw, I_h Tw) - a^1(Tw, Tw) \right)
    \end{align*}
    For this fixed $v$ we know that $a^1(I_h Tw, I_h Tw) -a^1(Tw,Tw) \rightarrow 0$ by the $H^1$ convergence assumption, and by $H^1$-continuity of $a^1$. This gives the result.
\end{proof}

\begin{theorem}\label{theorem:lowerbound}
    Assume that in addition to exactness on constant functions (i), $I_h$ satisfies both the $H^1$-convergence condition (ii) and the $L^2$-stability condition (iii), then
    \begin{align*}
        \liminf\limits_{h \rightarrow 0} \tilde \lambda_{h,1} \geq \lambda_1.
    \end{align*}
\end{theorem}

\begin{proof}
    Assume that $I_h$ satisfies both properties, but (for the sake of contradiction) that there exists
    a sequence $h_n \rightarrow 0$ monotonically with $\tilde \lambda_{h_n,1} \rightarrow \beta <  \lambda_1$.
    Note that by positivity of $\tilde a_h$ the sequence cannot diverge to $-\infty$,
    and by Theorem \ref{theorem:upperbound} it cannot diverge to $+\infty$.
    Hence, using ellipticity of $a^0$ on $\hat H^1$ we get that $|\tilde v_{h_n}|_{H^1} \leq C$ for
    some $C > 0$, otherwise the $a^0$ term would go to infinity. Here and below, $C$ represents a generic
    constant independent of $h_n$. By the Poincaré-Friedrichs theorem,
    it follows that
    \begin{align}\label{eq:bounded_evecs}
    \norm{\tilde v_{h_n}}_{H^1} \leq C
    \end{align}

    Recall that the Banach-Alaoglu theorem states that the closed unit ball in $H^1$ is weakly sequentially compact,
    while the Rellich-Kondrachev theorem asserts that norm-bounded $H^1$ sets are $L^2$-precompact.
    Hence we can assume that (going to a subsequence) $\tilde v_{h_n} \rightharpoonup \tilde v$ in $H^1$
    and $\tilde v_{h_n} \rightarrow \tilde v$ in $L^2$ for some $\tilde v \in H^1$.
    Note that the two limits are indeed the same since $\cdot \mapsto \Bra{\cdot,v}_{L^2}$
    is a continuous linear functional on $H^1$ for all $v \in L^2$.\\\\
    As $\tilde v_{h_n} \in \hat H^1$, it follows that $\tilde v \in \hat H^1$, the $L^2$-convergence gives us further that $\norm{v}_{L^2} = 1$. \\\\
    Without loss of generality, assume that $\ell^d(\Omega_1) = 1$ and define $J_n:\Omega_0\to \mathbb R$ as
    the constant function taking the value $\int_{\Omega_1} I_{h_n}T\tilde v_{h_n} d\ell^d$,
    and observe that
    \begin{align}\label{eq:image_bound}
    a^1(I_{h_n} T \tilde v_{h_n}, I_{h_n} T \tilde v_{h_n}) = a^1(I_{h_n} T \tilde v_{h_n} - J_n, I_{h_n}T \tilde v_{h_n} - J_n)
    \end{align}

    We know from the positivity of $a^0$ that $a^1(TI_{h_n}\tilde v_{h_n}, TI_{h_n}\tilde v_{h_n}) \leq \tilde r_{h_n}(\tilde v_{h_n})$
    and as the right hand side converges, the left hand side is bounded by a constant that does not depend on $n$.
    Hence, by ellipticity of $a^1$ on mean-free $H^1$ functions, it follows from \eqref{eq:image_bound} and the Poincaré-Friedrichs
    inequality that $\norm{I_{h_n} T \tilde v_{h_n} - J_n}_{H^1} < C$, yielding that $I_{h_n} T \tilde v_{h_n}$ is $H^1$-bounded
    as long as $J_n$ is -- which is the case by $L^2$-stability, the $L^2$-$L^2$ continuity of $T$, and the fact that $\norm{\tilde v_{h_n}}_{L^2} = 1$.
    Hence $I_{h_n} T \tilde v_{h_n}$ is (uniformly in $n$) $H^1$-bounded. Going to a subsequence,
    it is therefore weakly $H^1$- and strongly $L^2$-convergent (as above). As the weak-$H^1$ and $L^2$-limits coincide, to show
    that the sequence converges to $T\tilde v$ weakly in $H^1$, it is enough to show $L^2$-convergence.\\

    Let $C_1$ be the stability constant from the $L^2$-stability condition. We know that as $v_{h_n} \rightarrow \tilde v$ in $L^2$,
    it holds that $T\tilde v_{h_n} \rightarrow T \tilde v$. Therefore given $\varepsilon > 0$, pick $m$ so
    that $\norm{T\tilde v_{h_m} - T\tilde{v}}_{L^2} \leq \varepsilon/3$ and
    and $C_1 \norm{\tilde{v}_{h_n} - \tilde{v}_{h_m}}_{L^2} \leq \varepsilon/3$ for all $n \geq m$.
    Also, by the nesting property of the meshes $\tilde v_{h_n} - \tilde v_{h_m} \in \hat S_{h_n}^0$ and hence
    \begin{align*}
    \norm{I_{h_n} T \tilde v_{h_n} - T\tilde v}_{L^2}
    &\leq \norm{I_{h_n} T \tilde v_{h_n} - I_{h_n} T \tilde v_{h_m}}_{L^2} +
          \norm{I_{h_n} T \tilde v_{h_m} - T \tilde v_{h_m}}_{L^2} +
          \norm{T \tilde v_{h_m} - T \tilde v}_{L^2}\\
    &\leq C_1 \norm{\tilde v_{h_n} - \tilde v_{h_m}}_{L^2} +
              \norm{I_{h_n} T \tilde v_{h_m} - T \tilde v_{h_m} }_{L^2} +
              \norm{T \tilde v_{h_m} - T \tilde{v}}_{L^2} \\
    &\leq 2\varepsilon/3 + \norm{I_{h_n} T \tilde v_{h_m} - T \tilde v_{h_m}}_{L^2}.
    \end{align*}
    Now using the $H^1$-convergence property, pick $n_0 \geq m$ so that for $n \geq n_0$,
    we have $$\norm{I_{h_n} T \tilde v_{h_m} - T \tilde v_{h_m}}_{L^2} \leq \varepsilon/3.$$
    It follows that for all $n \geq n_0$
    \begin{align*}
    \norm{I_{h_n} T \tilde v_{h_n} - T \tilde v}_{L^2} \leq \epsilon.
    \end{align*}
    Summarizing, we now have proved that
    \begin{enumerate}
        \item $\tilde{v}_{h_n} \rightharpoonup \tilde{v}$ in $H^1$,
        \item $I_{h_n}T\tilde{v}_{h_n} \rightharpoonup T\tilde{v}$ in $H^1$, and
        \item $\tilde v \in \hat H^1$ with $\norm{\tilde{v}}_{L^2} = 1$
    \end{enumerate}
   The functions $u \mapsto \sqrt{a^0(u,u)}$ and $u \mapsto \sqrt{a^1(u,u)}$ are norms that are equivalent to the $\hat H^1$-norm by the Poincaré-Friedrichs inequality\footnote{The weak topologies also coincide, as for any $f \in \hat H^1$, the Riesz representation theorem gives a $\tilde f \in \hat H^1$ so that $\Bra{\cdot, f}_{H^1} = a_0(\cdot,\tilde f)$}.
   By the Banach-Steinhaus theorem, norms are weakly sequentially lower-semicontinuous. The $L^2$ convergence of $I_{h_n}Tv_{h_n}$ also gives that $J_n \rightarrow 0$ in $L^2$, and therefore also in $H^1$ as $J_n$ is constant).
   Thus we have:
\begin{align*}
    a^0(\tilde v,\tilde v) &\leq \liminf\limits_{n \rightarrow \infty}a^0(\tilde v_{h_n},\tilde v_{h_n}) \\
    a^1(T\tilde v,T\tilde v) &\leq \liminf\limits_{n \rightarrow \infty}a^1(I_{h_n}T\tilde v_{h_n} - J_n,I_{h_n}T\tilde v_{h_n} - J_n) \\
    &= \liminf\limits_{n \rightarrow \infty}a^1(I_{h_n} T \tilde v_{h_n}, I_{h_n}T\tilde v_{h_n})
\end{align*}
yielding
\begin{align*}
a(\tilde v,\tilde v)
    & = a^0(\tilde v,\tilde v) + a^1(T\tilde v,T\tilde v) \\
    & \leq \liminf\limits_{n \rightarrow \infty} \left( a^0(\tilde v_{h_n},\tilde v_{h_n}) + a^1(I_{h_n}T\tilde v_{h_n},I_{h_n}T\tilde v_{h_n}) \right) \\
    & = \liminf\limits_{n \rightarrow \infty} \tilde{r}_{h_n}(\tilde v_{h_n}) \\
    & = \lim\limits_{n \rightarrow \infty} \tilde \lambda_{h_n,1} = \beta < \lambda_1.
\end{align*}
    This is a contradiction to the defintion of $\lambda_1$, which requires that $a(\tilde{v},\tilde{v}) \geq \lambda_1$ and thus the claim is proved.
\end{proof}

\begin{theorem}
The conclusion of the previous theorem also holds for $I_h$ coming from the $H^1$-Galerkin approach.
\end{theorem}

\begin{proof}
The previous proof required the $L^2$-stability only in two places. The first was in bounding $J_n$, this is trivially bounded as $\norm{I_h}_{H^1\rightarrow H^1} = 1$ and $\tilde v_{h_n}$ is $H^1$-bounded. The second was in showing that $I_{h_n} T\tilde v_{h_n} \rightharpoonup T\tilde v$. But this follows from the fact that for $f \in H^1$:
\begin{align*}
\Bra{I_{h_n} T \tilde v_{h_n},f}_{H^1} = \Bra{T \tilde v_{h_n},I_{h_n} f}_{H^1}
\end{align*}
As $I_h f \rightarrow f$ in $H^1$ and $T \tilde v_{h_n} \rightharpoonup T \tilde v$, this gives the claim.
\end{proof}
\begin{theorem}[Eigenvector convergence]\label{theorem:evecs}
Assume that the eigenspace corresponding to $\lambda_1$ is one-dimensional. Then $|\Bra{\tilde v_h,v}_{L^2}| \rightarrow 1$ for $h \rightarrow 0$.
\end{theorem}
\begin{proof}
The proof of Theorem \ref{theorem:lowerbound} shows that for any monotone sequence $h_n \rightarrow 0$ for which $\lambda_{h_n}$ converges to $\lambda_1$,
every in $L^2$ convergent subsequence of $\tilde v_{h_n}$ converges to some $\tilde v \in H^1$ with $r(\tilde v) = \lambda_1$. Theorems \ref{theorem:lowerbound} and \ref{theorem:upperbound} show that $\lambda_{h_n} \rightarrow \lambda_1$ for any $h_n \rightarrow 0$. As $\tilde v_{h_n}$ is bounded in $H^1$ and therefore
pre-compact in $L^2$, this means that $\tilde v_{h_n} \rightarrow \tilde v$ in $L^2$ with $r(\tilde v) = \lambda_1$, from which it follows that $\tilde v = \pm v$, giving the result for a subsequence. We get the final result by noting that if $\tilde v_{h_n}$ has a subsequence $\tilde v_{h_{n_k}}$ for which $|\Bra{\tilde v_{h_{n_k}}, v}| \rightarrow c \neq 1$, by the above it must have a (further) subsequence that converges to a different value, a contradiction.
\end{proof}
\begin{theorem}[Higher eigenpairs]
The results of the previous theorems generalize to the $m$-th eigenpair for $m > 1$.
\end{theorem}
\begin{proof}
The generalization of Theorem \ref{theorem:upperbound} can be proven by the same method as the case $m=1$:
instead of chosing a $w$ with $r(w) < \lambda_1 + \epsilon$, choose a series of pairwise-orthogonal $(w_i)_{i=1\dots m} \in \hat S^0_h$ with the property that $r(w_i) \leq \lambda_i + \epsilon$.

The proof of Theorem \ref{theorem:lowerbound} must only be modified to show that the limit point $\tilde v_{m}$ is $L^2$-orthogonal to the previous eigenvectors $\spn\{v_1,\dots, v_{m-1}\}$. Writing $\tilde v_{m,h}$ for the $m$-th Ritz-vector of $\tilde a_h$ on $\hat S_h^1$, this readily follows (by induction) from the fact that $\tilde v_{m,h_n} \perp \spn\{\tilde v_{1,h_n},\dots,\tilde v_{m-1,h_n}\}$ and thus also $\Bra{\tilde v_{i,h_n},Tv_i} \rightarrow 0$ in $L^2$ for $i=1,\dots, m-1$.

The proof of the generalization of Theorem \ref{theorem:evecs} then follows from the generalization of Theorem \ref{theorem:lowerbound}.
\end{proof}
\paragraph{Convergence for homogeneous Dirichlet boundary.}
For the case of a homogeneous Dirichlet boundary, we must replace $\hat H^1$ with $H^1_0$.
Here we require the additional condition:
\begin{enumerate}[(i')]
\item Boundary preservation: for $v \in H^1_0(\Omega_1) \cap S_h^0$ it holds that  $I_h Tv \in H^1_0(\Omega_1)$
\end{enumerate}
We note that the $I_h$ from collocation and $H^1$-Galerkin have this property.
\begin{theorem}
The same convergence results hold also for the Ritz-values for homogeneous Dirichlet boundary conditions, provided that $I_h$ satisfies (i') in addition
to satisfying (ii) and (iii) from the previous section.
\end{theorem}
\begin{proof}
Same as the proofs for the Neumann case, except that coercivity of the bilinear forms involved on $\hat H^1$ is replaced with that on $H^1_0$. The boundary preservation property ensures we do not need $J_n$ in the proof of Theorem \ref{theorem:lowerbound}, the proof otherwise goes on exactly the same lines.
\end{proof}

\subsection{Numerical Experiments}
The three TO approaches yield a stiffness-matrix $D$ of the form
\begin{align*}
D &= D_0 + A^TD_1A
\end{align*}
where $A$ is the representation matrix of $I_hT$ and $D_0,D_1$ result from the bilinear form of the static Laplacian.
In the case of the $L^2$- or $H^1$-Galerkin approach, the matrix $A$ has the form
\begin{align*}
A &= G^{-1}\tilde A,
\end{align*}
where $\tilde A = \Bra{\varphi_i, T\varphi_j}$ and $G_{i,j} = \Bra{\varphi_i,\varphi_j}$.
Note that sparsity of $G$ and $\tilde A$ does not necessarily imply sparsity of $A$. A naive calculation of the full matrix $A$ therefore renders the computation of (some) eigenvalues of $D$ too expensive for larger problems.  We therefore do not include the Galerkin-TO approaches in the numerical experiments below, further work is needed to determine if the Galerkin-TO methods can be modified to overcome this issue.

All TO results shown below are therefore the collocation-based ones.
For the ``non-adaptive TO'' experiments, we used identical regular meshes for the initial and final triangulations $\mathcal T_h^0$ and $\mathcal T_h^1$.
For the ``adaptive TO'' experiments, the initial mesh was regular but the final mesh was a Delaunay-triangulation
(with the \texttt{VoronoiDelaunay.jl} package \cite{voronoidelaunay})
of the images of the initial mesh points under the flow map $T$.
The theory outlined above requires (i) nested meshes, and (ii) uniform bounds on shape-regularity and quasi-uniformity constants of the meshes as they are refined;  each of these properties is difficult to guarantee in general for the image triangulations.
We nevertheless include the adaptive TO in the plots for comparison.

We note that the finest grid used for the TO experiments is $257\times257$.
This is not as fine as the finest grid used in the CG approach ($513\times 513)$.
The reason for this is that the stiffness matrix $D$ from the TO approach does not have the same banded structure as that coming from the CG approach,
making a solution of the eigenproblem more expensive (for the example of the cylinder flow: 1h 28min for the non-adaptive TO using a $P^1$ $513\times513$ mesh vs.\ 14 minutes for the CG approach with $P^2$ elements on a $1025\times1025$ mesh).  The reference solution here is the same as that used for the CG experiments.

\subsubsection{Standard Map}

We investigate the convergence of collocation-based TO approaches on the Standard map test case.
\begin{figure}[H]
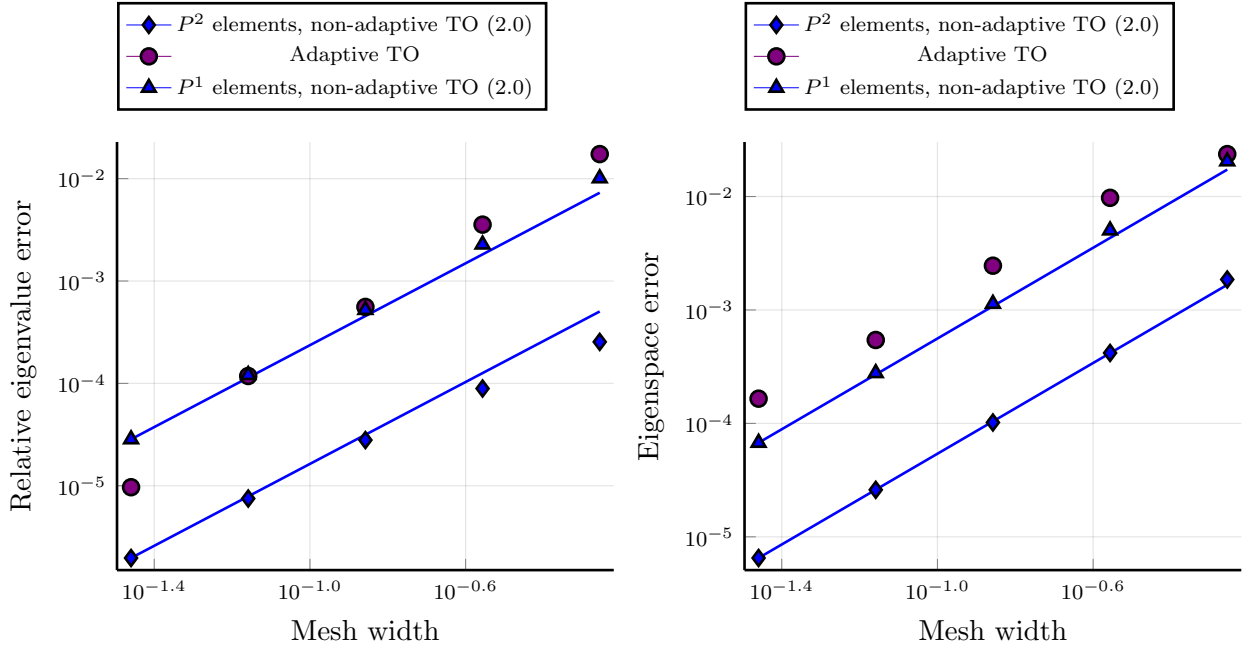

\begin{subfigure}{0.5\linewidth}
\input{Figures/processed/SM2_eigval_errs2_TO.tex}
\end{subfigure}
\begin{subfigure}{0.5\linewidth}
\input{Figures/processed/SM2evec2errsTO.tex}
\end{subfigure}
\caption{Standard map: Errors in the first nontrivial eigenvalue (left) and  2-dimensional eigenspace spanned by the smallest eigenvectors (right) of the dynamic Laplacian discretized with the TO approach. The slopes of the corresponding lines are given in brackets in the legends.}
\label{fig:std_map_TO}
\end{figure}
We see an improvement by a full order of magnitude in the errors in both the eigenvalues and the eigenspaces for $P^2$ elements over $P^1$ elements.
Concerning the convergence order, this experiment suggests that the order for $P^2$ elements is not higher than for $P^1$ elements in the TO approach, which match those obtained for eigenvectors based on $P^1$ elements in the CG approach.

\subsubsection{Cylinder Flow}

\begin{figure}[H]
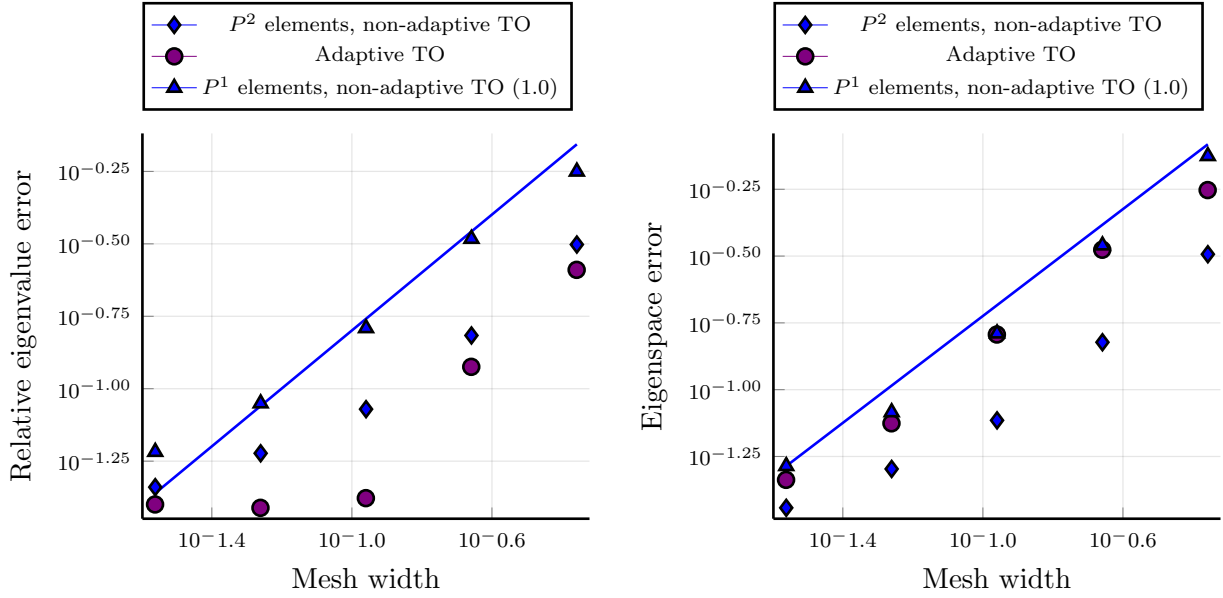

\begin{subfigure}{0.45\linewidth}
\input{Figures/processed/CYq58_eigval_errs2_TO.tex}
\end{subfigure}
\qquad
\begin{subfigure}{0.45\linewidth}
\input{Figures/processed/CYq58evec2errsTO.tex}
\end{subfigure}
\caption{Cylinder flow: Errors in the first nontrivial eigenvalue (left) and corresponding eigenspace (right) of the dynamic Laplacian for the TO approaches. The slopes of the corresponding lines are given in brackets in the legends.}
\label{fig:rot_dg_10_TO}
\end{figure}

We observe a smaller error reduction in the eigenvalue and eigenvector errors for $P^2$ vs $P^1$ elements for the cylinder flow (Fig.~\ref{fig:rot_dg_10_TO}) compared to the standard map example.  The reduction in convergence orders compared to the standard map (i.e.\ comparing Fig.~\ref{fig:rot_dg_10_TO} to Fig.~\ref{fig:std_map_TO}) mirrors the reductions observed in the corresponding CG experiments (i.e.\ comparing Fig.~\ref{fig:cyl_CG} to Fig.~\ref{fig:std_map_CG}), where the orders are around half those seen for the standard map.
Similarly to the standard map results, the $P^2$ elements did not perform asymptotically better than $P^1$ elements.

\subsubsection{Bickley Jet}

\begin{figure}[H]
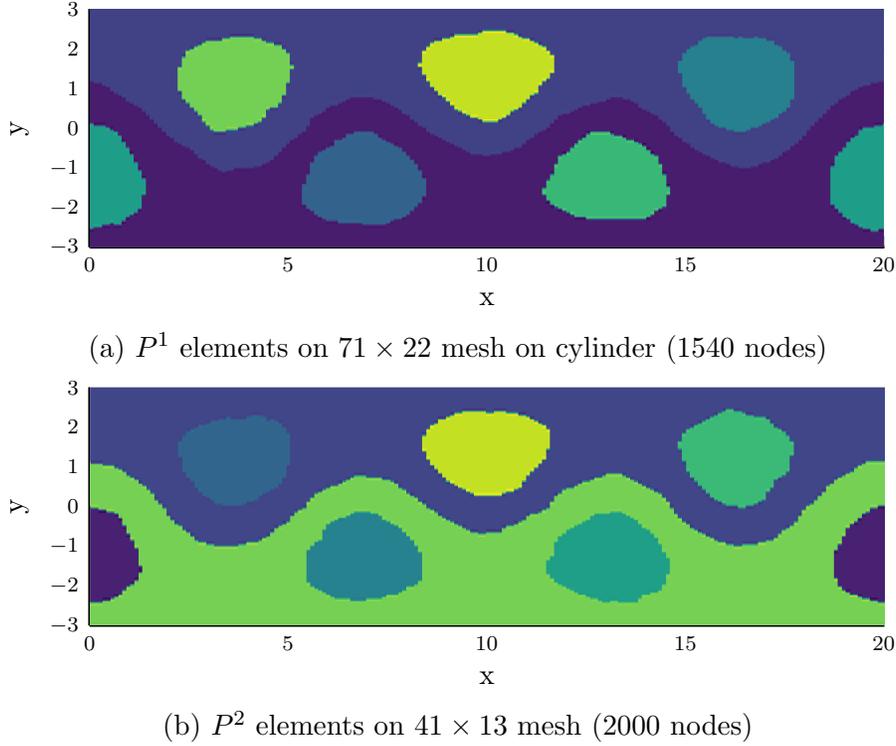

\centering
\begin{subfigure}{0.8\linewidth}
\centering
\input{Figures/processed/bickleyaxestop1.tex}
\subcaption{$P^1$ elements on $71\times 22$ mesh on cylinder (1540 nodes)}
\end{subfigure}
\begin{subfigure}{0.8\linewidth}
\centering
\input{Figures/processed/bickleyaxestop2.tex}
\subcaption{$P^2$ elements on $41\times 13$ mesh (2000 nodes)}
\end{subfigure}
\caption{Bickley jet: Comparison of TO method with $P^1$ and $P^2$ elements.}
\label{fig:bickley_to}
\end{figure}

We repeat the experiment from Section~\ref{eq:exp_CG_Bickley} with the (non-adaptive) TO approach, see \ Fig.~\ref{fig:bickley_to}. 
As in Section~\ref{eq:exp_CG_Bickley}, the mesh widths in Fig.~\ref{fig:bickley_to} below are the lowest for which the topology of the clustering result does not change.  Here, the dominant computational cost is given by the evaluation of the flow map, which has to be evaluated once for each basis function of the finite element space. Note that there are three basis functions per element for $P^1$ and six per element for $P^2$ elements.  While for $P^2$ elements, a coarser mesh was sufficient, the number of basis functions was comparable (even a little larger) than for the $P^1$ case -- which is in contrast to the corresponding CG experiment.

\subsection{TO Convergence Rates}

Determining the theoretical convergence orders of the TO-methods is an outstanding task.  In order to investigate whether we can expect higher-order convergence, we consider the collocation-based non-adaptive TO approach on a simple example in one dimension.
In this example we can largely isolate the dynamics from the errors to focus on errors arising from translations of the basis functions.

The only volume-preserving diffeomorphisms of the circle $\mathbb S^1 = \mathbb R / \mathbb Z$ to itself are rigid rotations.
Let $T : \mathbb S^1 \rightarrow \mathbb S^1$ be given by $T(x) = (x + \alpha) \mod 1$.
Rigid rotations commute with the Laplace operator;  thus $\Delta^{dyn}=(\Delta+T^*_1\Delta T_{1,*})/2=\Delta$,
and the dynamic Laplacian is equal to the static Laplacian.
The first nontrivial eigenspace of the dynamic Laplacian is therefore spanned by $v(x) = \sin(2\pi x)$ and $u(x) = \cos(2\pi x)$ which are orthogonal eigenvectors for the eigenvalue $\lambda_1 = 4\pi^2$.

\begin{figure}[H]
\begin{subfigure}{0.5\linewidth}

\begin{tikzpicture}[]
\begin{axis}[scale only axis,width=0.8\linewidth, legend pos = {south east}, ylabel = {Relative eigenvalue error}, xmin = {0.0017602548097867773}, xmax = {0.0693480920042403}, ymax = {0.0018679494904262305}, ymode = {log}, xlabel = {Mesh width}, unbounded coords=jump,scaled x ticks = false,xlabel style = {font = {\fontsize{11 pt}{14.3 pt}\selectfont}, color = {rgb,1:red,0.00000000;green,0.00000000;blue,0.00000000}, draw opacity = 1.0, rotate = 0.0},log basis x=2,xmajorgrids = true,xtick = {0.001953125,0.00390625,0.0078125,0.015625,0.03125,0.0625},xticklabels = {$2^{-9}$,$2^{-8}$,$2^{-7}$,$2^{-6}$,$2^{-5}$,$2^{-4}$},xtick align = inside,xticklabel style = {font = {\fontsize{8 pt}{10.4 pt}\selectfont}, color = {rgb,1:red,0.00000000;green,0.00000000;blue,0.00000000}, draw opacity = 1.0, rotate = 0.0},x grid style = {color = {rgb,1:red,0.00000000;green,0.00000000;blue,0.00000000},
draw opacity = 0.1,
line width = 0.5,
solid},axis lines* = left,x axis line style = {color = {rgb,1:red,0.00000000;green,0.00000000;blue,0.00000000},
draw opacity = 1.0,
line width = 1,
solid},scaled y ticks = false,ylabel style = {font = {\fontsize{11 pt}{14.3 pt}\selectfont}, color = {rgb,1:red,0.00000000;green,0.00000000;blue,0.00000000}, draw opacity = 1.0, rotate = 0.0},log basis y=2,ymajorgrids = true,ytick = {4.76837158203125e-7,3.814697265625e-6,3.0517578125e-5,0.000244140625},yticklabels = {$2^{-21}$,$2^{-18}$,$2^{-15}$,$2^{-12}$},ytick align = inside,yticklabel style = {font = {\fontsize{8 pt}{10.4 pt}\selectfont}, color = {rgb,1:red,0.00000000;green,0.00000000;blue,0.00000000}, draw opacity = 1.0, rotate = 0.0},y grid style = {color = {rgb,1:red,0.00000000;green,0.00000000;blue,0.00000000},
draw opacity = 0.1,
line width = 0.5,
solid},axis lines* = left,y axis line style = {color = {rgb,1:red,0.00000000;green,0.00000000;blue,0.00000000},
draw opacity = 1.0,
line width = 1,
solid},    xshift = 0.0mm,
    yshift = 0.0mm,
    axis background/.style={fill={rgb,1:red,1.00000000;green,1.00000000;blue,1.00000000}}
,legend style = {color = {rgb,1:red,0.00000000;green,0.00000000;blue,0.00000000},at={(0.0,1.2)},anchor=west,
draw opacity = 1.0,
line width = 1,
solid,fill = {rgb,1:red,1.00000000;green,1.00000000;blue,1.00000000},font = {\fontsize{8 pt}{10.4 pt}\selectfont}},colorbar style={title=}, xmode = {log}, ymin = {8.91294606100333e-8}]\addplot+[draw=none, color = {rgb,1:red,0.00000000;green,0.00000000;blue,1.00000000},
draw opacity = 1.0,
line width = 0,
solid,mark = triangle*,
mark size = 3.0,
mark options = {
    color = {rgb,1:red,0.00000000;green,0.00000000;blue,0.00000000}, draw opacity = 1.0,
    fill = {rgb,1:red,0.00000000;green,0.00000000;blue,1.00000000}, fill opacity = 1.0,
    line width = 1,
    rotate = 0,
    solid
}] coordinates {
(0.0625, 0.0005794621986196586)
(0.03125, 0.0014094930677278873)
};
\addlegendentry{$P^1$ elements non-adaptive TO (2.0)}
\addplot+[draw=none, color = {rgb,1:red,0.00000000;green,0.00000000;blue,1.00000000},
draw opacity = 1.0,
line width = 0,
solid,mark = triangle*,
mark size = 3.0,
mark options = {
    color = {rgb,1:red,0.00000000;green,0.00000000;blue,0.00000000}, draw opacity = 1.0,
    fill = {rgb,1:red,0.00000000;green,0.00000000;blue,1.00000000}, fill opacity = 1.0,
    line width = 1,
    rotate = 0,
    solid
},forget plot] coordinates {
(0.0078125, 8.833482455916898e-5)
(0.00390625, 2.0089840120034747e-6)
(0.001953125, 5.521871608029213e-6)
};
\addplot+[draw=none, color = {rgb,1:red,0.00000000;green,0.00000000;blue,1.00000000},
draw opacity = 1.0,
line width = 0,
solid,mark = diamond*,
mark size = 3.0,
mark options = {
    color = {rgb,1:red,0.00000000;green,0.00000000;blue,0.00000000}, draw opacity = 1.0,
    fill = {rgb,1:red,0.00000000;green,0.00000000;blue,1.00000000}, fill opacity = 1.0,
    line width = 1,
    rotate = 0,
    solid
}] coordinates {
(0.0625, 0.00015547774184128594)
(0.03125, 3.0166815377862738e-5)
(0.015625, 9.874041872725751e-6)
(0.0078125, 1.8895805447368266e-6)
(0.00390625, 6.177948000977439e-7)
(0.001953125, 1.1812000664668654e-7)
};
\addlegendentry{$P^2$ elements non-adaptive TO}
\addplot+ [color = {rgb,1:red,0.00000000;green,0.00000000;blue,1.00000000},
draw opacity = 1.0,
line width = 1,
solid,mark = none,
mark size = 2.0,
mark options = {
    color = {rgb,1:red,0.00000000;green,0.00000000;blue,0.00000000}, draw opacity = 1.0,
    fill = {rgb,1:red,0.00000000;green,0.00000000;blue,1.00000000}, fill opacity = 1.0,
    line width = 1,
    rotate = 0,
    solid
},forget plot]coordinates {
(0.0625, 0.00015547774184128602)
(0.03125, 3.8869435460321505e-5)
(0.015625, 9.717358865080376e-6)
(0.0078125, 2.429339716270094e-6)
(0.00390625, 6.073349290675235e-7)
(0.001953125, 1.5183373226688088e-7)
};
\end{axis}

\end{tikzpicture}
\end{subfigure}
\begin{subfigure}{0.5\linewidth}
\begin{tikzpicture}[]
\begin{axis}[scale only axis,width=0.8\linewidth,legend pos = {south east}, ylabel = {Eigenvector error}, xmin = {0.0017602548097867773}, xmax = {0.0693480920042403}, ymax = {0.007687854899346467}, ymode = {log}, xlabel = {Mesh width}, unbounded coords=jump,scaled x ticks = false,xlabel style = {font = {\fontsize{11 pt}{14.3 pt}\selectfont}, color = {rgb,1:red,0.00000000;green,0.00000000;blue,0.00000000}, draw opacity = 1.0, rotate = 0.0},log basis x=2,xmajorgrids = true,xtick = {0.001953125,0.00390625,0.0078125,0.015625,0.03125,0.0625},xticklabels = {$2^{-9}$,$2^{-8}$,$2^{-7}$,$2^{-6}$,$2^{-5}$,$2^{-4}$},xtick align = inside,xticklabel style = {font = {\fontsize{8 pt}{10.4 pt}\selectfont}, color = {rgb,1:red,0.00000000;green,0.00000000;blue,0.00000000}, draw opacity = 1.0, rotate = 0.0},x grid style = {color = {rgb,1:red,0.00000000;green,0.00000000;blue,0.00000000},
draw opacity = 0.1,
line width = 0.5,
solid},axis lines* = left,x axis line style = {color = {rgb,1:red,0.00000000;green,0.00000000;blue,0.00000000},
draw opacity = 1.0,
line width = 1,
solid},scaled y ticks = false,ylabel style = {font = {\fontsize{11 pt}{14.3 pt}\selectfont}, color = {rgb,1:red,0.00000000;green,0.00000000;blue,0.00000000}, draw opacity = 1.0, rotate = 0.0},log basis y=2,ymajorgrids = true,ytick = {3.814697265625e-6,3.0517578125e-5,0.000244140625,0.001953125},yticklabels = {$2^{-18}$,$2^{-15}$,$2^{-12}$,$2^{-9}$},ytick align = inside,yticklabel style = {font = {\fontsize{8 pt}{10.4 pt}\selectfont}, color = {rgb,1:red,0.00000000;green,0.00000000;blue,0.00000000}, draw opacity = 1.0, rotate = 0.0},y grid style = {color = {rgb,1:red,0.00000000;green,0.00000000;blue,0.00000000},
draw opacity = 0.1,
line width = 0.5,
solid},axis lines* = left,y axis line style = {color = {rgb,1:red,0.00000000;green,0.00000000;blue,0.00000000},
draw opacity = 1.0,
line width = 1,
solid},    xshift = 0.0mm,
    yshift = 0.0mm,
    axis background/.style={fill={rgb,1:red,1.00000000;green,1.00000000;blue,1.00000000}}
,legend style = {color = {rgb,1:red,0.00000000;green,0.00000000;blue,0.00000000},at={(0.0,1.2)},anchor=west,
draw opacity = 1.0,
line width = 1,
solid,fill = {rgb,1:red,1.00000000;green,1.00000000;blue,1.00000000},font = {\fontsize{8 pt}{10.4 pt}\selectfont}},colorbar style={title=}, xmode = {log}, ymin = {5.0290217043355e-7}]\addplot+[draw=none, color = {rgb,1:red,0.00000000;green,0.00000000;blue,1.00000000},
draw opacity = 1.0,
line width = 0,
solid,mark = triangle*,
mark size = 3.0,
mark options = {
    color = {rgb,1:red,0.00000000;green,0.00000000;blue,0.00000000}, draw opacity = 1.0,
    fill = {rgb,1:red,0.00000000;green,0.00000000;blue,1.00000000}, fill opacity = 1.0,
    line width = 1,
    rotate = 0,
    solid
}] coordinates {
(0.0625, 0.0058530330942493216)
(0.03125, 0.0014433907916510517)
};
\addlegendentry{$P^1$ elements non-adaptive TO (2.0)}
\addplot+[draw=none, color = {rgb,1:red,0.00000000;green,0.00000000;blue,1.00000000},
draw opacity = 1.0,
line width = 0,
solid,mark = triangle*,
mark size = 3.0,
mark options = {
    color = {rgb,1:red,0.00000000;green,0.00000000;blue,0.00000000}, draw opacity = 1.0,
    fill = {rgb,1:red,0.00000000;green,0.00000000;blue,1.00000000}, fill opacity = 1.0,
    line width = 1,
    rotate = 0,
    solid
},forget plot] coordinates {
(0.0078125, 8.982518446045041e-5)
(0.00390625, 2.2451468429275467e-5)
(0.001953125, 5.612573493039485e-6)
};
\addplot+[draw=none, color = {rgb,1:red,0.00000000;green,0.00000000;blue,1.00000000},
draw opacity = 1.0,
line width = 0,
solid,mark = diamond*,
mark size = 3.0,
mark options = {
    color = {rgb,1:red,0.00000000;green,0.00000000;blue,0.00000000}, draw opacity = 1.0,
    fill = {rgb,1:red,0.00000000;green,0.00000000;blue,1.00000000}, fill opacity = 1.0,
    line width = 1,
    rotate = 0,
    solid
}] coordinates {
(0.0625, 0.0009321517460575359)
(0.03125, 0.00018086095442021853)
(0.015625, 5.5459349261077855e-5)
(0.0078125, 1.0612520522592156e-5)
(0.00390625, 3.4546154743775234e-6)
(0.001953125, 6.605530590042567e-7)
};
\addlegendentry{$P^2$ elements non-adaptive TO (2.0)}
\addplot+ [color = {rgb,1:red,0.00000000;green,0.00000000;blue,1.00000000},
draw opacity = 1.0,
line width = 1,
solid,mark = none,
mark size = 2.0,
mark options = {
    color = {rgb,1:red,0.00000000;green,0.00000000;blue,0.00000000}, draw opacity = 1.0,
    fill = {rgb,1:red,0.00000000;green,0.00000000;blue,1.00000000}, fill opacity = 1.0,
    line width = 1,
    rotate = 0,
    solid
},forget plot]coordinates {
(0.0625, 0.005853033094249322)
(0.03125, 0.0014632582735623315)
(0.015625, 0.00036581456839058287)
(0.0078125, 9.145364209764572e-5)
(0.00390625, 2.286341052441143e-5)
(0.001953125, 5.715852631102857e-6)
};
\addplot+ [color = {rgb,1:red,0.00000000;green,0.00000000;blue,1.00000000},
draw opacity = 1.0,
line width = 1,
solid,mark = none,
mark size = 2.0,
mark options = {
    color = {rgb,1:red,0.00000000;green,0.00000000;blue,0.00000000}, draw opacity = 1.0,
    fill = {rgb,1:red,0.00000000;green,0.00000000;blue,1.00000000}, fill opacity = 1.0,
    line width = 1,
    rotate = 0,
    solid
},forget plot]coordinates {
(0.0625, 0.0009321517460575358)
(0.03125, 0.00023303793651438396)
(0.015625, 5.825948412859599e-5)
(0.0078125, 1.4564871032148997e-5)
(0.00390625, 3.6412177580372493e-6)
(0.001953125, 9.103044395093123e-7)
};
\end{axis}

\end{tikzpicture}
\end{subfigure}
\caption{1D shift map: Errors in the first calculated eigenvalue (left) and corresponding eigenvector $\tilde v_{1,h}$. For the case $h=2^{-6}$, the call to \texttt{eigs} failed for the $P^1$ non-adaptive TO.}
\label{fig:shift_map_rates}
\end{figure}
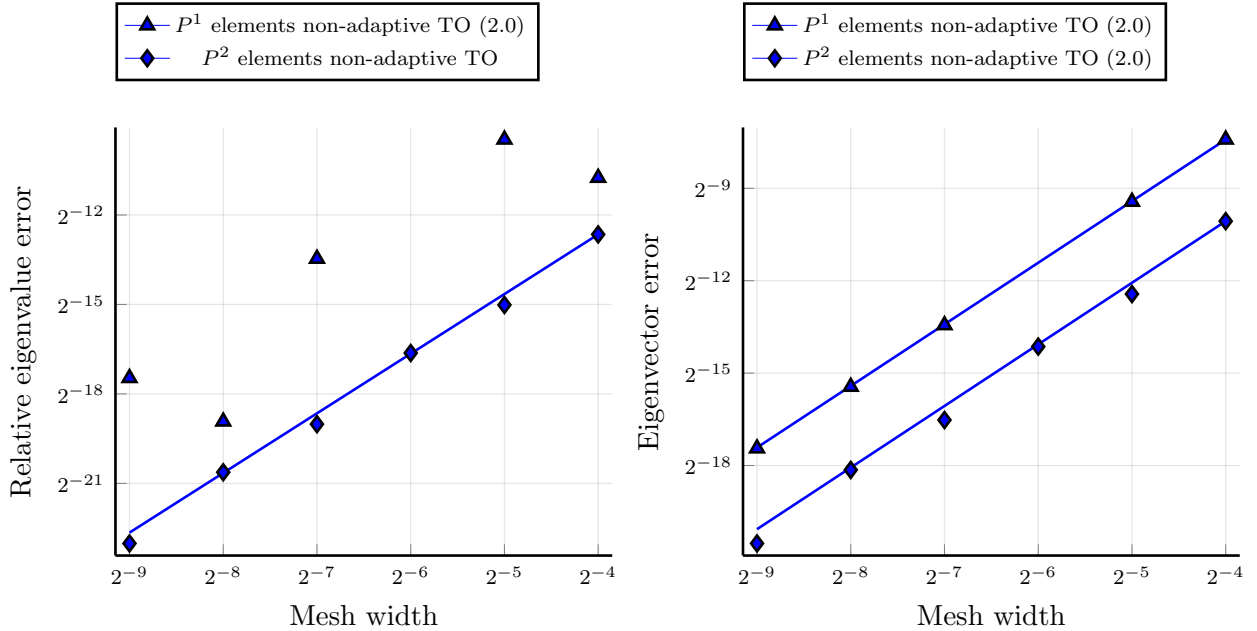

Any errors in the corresponding discrete bilinear form
$$\tilde a_h(u,v) = \frac12(a(u,v) + a(I_h Tu, I_h Tv))$$
therefore arise solely from discretisation errors related to the rotation of the $\varphi_i$ by $\alpha$.
In our experiments, we consider $\alpha = 0.15$.
For the non-adaptive collocation TO approach, denote the leading nontrivial eigenvector of $\Delta^{dyn}$ by $\tilde v_{h,1}$.
We now look at the order of convergence for the error $\|\tilde v_{h,1}-v\|_{L^2}$ as $h \rightarrow 0$. 
 We compute the errors independently to the numerical experiments by using Fourier coefficients.  
 With the Fourier coefficients $A_h = \Bra{v, \tilde v_{1,h}}_{L^2}$ and $B_h = \Bra{u, \tilde v_{1,h}}_{L^2}$ of $\tilde v_{h,1}$, $P v_{h,1} = A_h v + B_h u$
is the $L^2$-orthogonal projection of $\tilde v_{1,h}$ onto the first nontrivial eigenspace. As a measure of eigenvector error we now choose $e_h = \|P \tilde  v_{h,1} - \tilde v_{1,h}\|_{L^2}^2$.
All integrals were approximated by the trapezoid rule with $10^8$ quadrature points on $[0,1]$;  see Figure \ref{fig:shift_map_rates}.

The convergence rates in Figure \ref{fig:shift_map_rates} are identical to those found for the standard map in Figure \ref{fig:std_map_CG} and Figure \ref{fig:std_map_TO}. 
The CG approach in this example is simply computing the eigenvectors and eigenvectors of the standard Laplace operator and therefore unsurprisingly one recovers the theoretical orders of convergence.
These numerical results suggest that when using the non-adaptive TO method, we cannot
expect an asymptotically higher convergence rate for $P^2$ elements in comparison to $P^1$-Elements even for very simple flow maps.

\section{Conclusion}
We compared the use of $P^1$ and $P^2$ elements in collocation-based CG and TO approximations of the
dynamic Laplacian. 
In the CG approach applied to weakly nonlinear dynamics, $P^2$ elements can significantly reduce the computational cost by providing an asymptotically higher order of convergence.
A benefit of $P^1$ elements for the CG approach
is that a first-order method of quadrature can be used in the discretization, whereas using first-order quadrature for $P^2$ elements results in a singular mass matrix.
This does not affect the asymptotic order of convergence, but nevertheless introduces a constant factor which may be relevant in some low-data cases.

In the non-adaptive TO approach, there seems to be little benefit gained by using $P^2$ elements as opposed to $P^1$ elements for the non-adaptive TO approximation. 
The adaptive TO is inherently $P^1$-based, and thus
does not benefit from a $P^2$ discretization either.
In general, numerical experiments suggest that collocation-based TO approaches have equal (or in some cases better) rates of convergence compared to the $P^1$ CG discretisation. 

The $P^2$ CG discretisation generally had a higher order of convergence, though the numerically observed convergence rates varied significantly when the
quadrature order was changed or when applied to more nonlinear and hyperbolic dynamics. 
It is difficult to compare the CG and TO approaches in general as the former relies on being able to calculate derivatives of the flow map, whereas the latter is purely data based.
This makes
the TO approach applicable to some cases where the CG approach cannot be used.

A hindrance to using the TO approach for finely resolved meshes is the fact that here the call to \texttt{eigs} takes much longer
compared to when one uses the CG discretisation on the same mesh. 
We suspect that this is due to the fact that unlike in the TO approach, the CG discretisation
preserves the banded structure of the stiffness matrix. More work is needed to determine how the eigenproblem can be solved efficiently in this case, or whether it is possible to avoid the eigenproblem completely but still be able to compute coherent sets. More work is also needed to determine the true rates of convergence of eigenvalues and eigenvectors. We proved that they do converge for the non-adaptive TO approach, but were only
able to conjecture what the true rates are based on numerical experiments. It also remains to be seen whether Galerkin TO approaches can be modified to be computationally efficient.

\section{Acknowledgements}
OJ and NS acknowledge support by the Priority Programme SPP 1881 Turbulent Superstructures of the Deutsche Forschungsgemeinschaft and a Universities Australia / DAAD travel grant.
GF is partially supported by an ARC Discovery Project and a joint Universities Australia / DAAD travel award.
GF thanks the Faculty of Mathematics at the Technical University Munich for hospitality during his research visits.
We would like thank Christian Ludwig for the reference to \cite{cai95} used in the proof of Lemma \ref{lemma:stability}.
We also thank Alavaro de Diego and Daniel Karrasch for helpful comments and their contributions to the \texttt{CoherentStructures.jl} package.
\begin{appendix}

\section{Proofs}
Throughout all proofs, $C$ refers to a constant depending only on the mesh and dimension.

\subsection{Proof of Lemma \ref{lemma:stability}}

The following proof of Lemma \ref{lemma:stability} uses ideas from the proof of Lemma 2.1 from \cite{cai95}.
We start with a helper lemma:
\begin{lemma}\label{lemma:overlappingtriangles}
Let $\cT^0_{h'}$ and $\cT^1_{h}$ be quasi-uniform meshes on $\Omega_0$ and $\Omega_1$ respectively for some values of $h,h' > 0$. Let $F : \Omega_0 \rightarrow \Omega_1$ be smooth on $\Omega_0$ with smooth extension to the boundary. Then there exists $C > 0$ independent of $h$ and $h'$ so that for all $\tau' \in \cT^0_{h'}$ it holds that:
\begin{align*}
\mathcal N(\tau') := |\{ \tau \in \cT_{h}^1 : F(\tau') \cap \tau \neq \emptyset\}| \leq C\left( \frac{h' + h}{h}\right)^d
\end{align*}
\end{lemma}
\begin{proof}
By quasi-uniformity of the meshes, we can find positive constants $a',A'$ so that for any triangle $\tau' \in \cT^0_{h'}$ there is a point $\tau'_x$ so that:
\begin{align}\label{eq:ballcondition}
B(\tau'_x, a'h') \subset \tau \subset B(\tau'_x, A'h')
\end{align}
This gives that where $y := F(\tau'_x)$:
\begin{align}
F(\tau') \subset F(B(\tau'_x, A'h')) \subset B(y, \underbrace{\norm{DF}_\infty A'}_{:=E}h' )
\end{align}
where the last inclusion follows from the mean-value theorem,
and $\norm{DF}_\infty$ exists as $F$ is smoothly extensible to the boundary.\\

Now write $a,A$ for the
shape-regularity/quasi-uniformity constants of $\cT_h^1$,
so that a formula like \eqref{eq:ballcondition} holds
for $\cT_h^1$ also.\\

Then if $\tau \in \mathcal T_h^1$ intersects $B(y, E h')$, it must hold that $B(\tau_x, ah) \subset B(\tau_x, Ah) \subset B(y, Eh' + 2Ah) \subset B(y, c(h' + h)) $ with suitable $c > 0$. As for different $\tau$, the sets $B(\tau_x, ah)$ are disjoint, we get
\begin{align*}
\mathcal N (\tau') \leq \left(\frac{c(h' + h)}{ah}\right)^d
\end{align*}
which gives the claim.
\end{proof}

\begin{lemma}\label{lemma:massmatrixspectrum}
Let $\cT_h$, $h>0$ be a family of quasi-uniform meshes on an open subset of $d$-dimensional space. Let $S_h$ be the space of functions representable by $P^k$-Lagrange elements on the mesh. Then there exist $C,C' > 0$ so that for all $v \in S^h$:
\begin{align*}
C_1 \norm{v}_{L^2}^2 \leq h^d \sum_p |v(p)|^2 \leq C_2 \norm{v}_{L^2}^2
\end{align*}
where $p$ ranges over the nodes of the triangulation.
\end{lemma}
\begin{proof}
This follows directly from well-known results about the spectrum of mass-matrices, see \cite[p.386]{ernguermond}
\end{proof}
\begin{proof}[Proof of Lemma \ref{lemma:stability}]
Throughout the proof, $C$ refers to a constant that does
not depend on $h$ or $h'$ and whose exact value can change from line to line.
Let $h' \geq h > 0$ and $v \in S^0_{h'}$.
Then by Lemma \ref{lemma:massmatrixspectrum}
\begin{align*}
\norm{I_h Tv}_{L^2}^2 &\leq Ch^d \sum_p |v(F^{-1}(p))|^2 \eqqcolon (*)
\end{align*}
where $p$ ranges over the nodes of $\cT^1_{h}$.
Using Lemma \ref{lemma:overlappingtriangles},
we know that for any triangle $\tau' \in \mathcal T^0_h$,
at most $C\left(\frac{h'+h}{h}\right)^d$ triangles from $\mathcal T^1_{h'}$ can intersect with $F(\tau')$, up to a constant factor this therefore bounds the number of vertices $p$ in $\mathcal T^1_h$ for which $F^{-1}(p)$ lies in a given triangle.
Moreover, as we are using $P^1$-Lagrange elements, $|v(F^{-1}(p))|$ is bounded
by $|v(p')|$ for some vertex $p'$ of the triangle that contains $F^{-1}(p)$. This gives, using Lemma \ref{lemma:massmatrixspectrum} and the fact that $h' \geq h$ that:
\begin{align*}
        (*)            &\leq C h^d \left(\frac{h' + h}{h}\right)^d \sum_{p'} |v(p')|^2 \\
                    &\leq C (h')^d  \sum_{p'} |v(p')|^2 \\ &\leq C \norm{v}_{L^2}^2
\end{align*}
where $p'$ ranges over the nodes of $\cT^0_{h'}$.
\end{proof}

\subsection{Proof of Lemma \ref{lemma:convergence}}

\begin{proof}
Without loss of generality, look at triangles in dimension $2$. Since $\norm{I_h v -v}_{L^2} \leq C \norm{I_h v - v}_{L^\infty}$ and $v$ is uniformly continuous, we have that
    \[
    \norm{I_h v - v}_{L^2} \rightarrow 0 \qquad \text{as } h\to 0.
    \]

    What remains to be shown is $\norm{\nabla(I_h v -v)}_{L^2} \rightarrow 0$.

    We first show that $\norm{\nabla I_h v}_{L^\infty} \leq C \norm{\nabla v}_{L^\infty}$ (where $C$ does not depend on $h$).
    It is enough to prove this for any triangle $\tau$ in the mesh $\cT^1_h$.
    Assume first that the triangle $\tau$ has vertices $0, e_1$ and $e_2$, without loss of generality also $v(0) = 0$.
    Then $\nabla I_h v  = (v(e_1), v(e_2))^T$.

 By a mean value inequality\footnote{As $v$ is only piecewise $C^\infty$, we cannot use the mean value theorem directly. However, as we only care about a weaker statement of the form $|f(1) - f(0)| \leq \norm{f'}_{L^\infty([0,1])}$  this is not a problem.
 },
 $\max\{|v(e_1)|,|v(e_2)|\} \leq \norm{\nabla v}_{L^\infty(\tau)}$.
 Thus $\norm{\nabla I_h v}_{L^\infty(\tau)} \leq \norm{\nabla v}_{L^\infty(\tau)}$. Shape regularity immediately gives $\norm{\nabla I_h v}_{L^\infty(\tau)} \leq C\norm{\nabla v}_{L^\infty(\tau)}$ for general triangles, taking suprema over all triangles gives the claim.

    By assumption $v$ is $C^2$ except for on a nowhere dense set of measure zero. By standard FEM theory, $\norm{\nabla (I_h v - v)}_{L^2(\Omega')} \rightarrow 0$ on all sub-meshes $\Omega'$ on which $v$ is $C^\infty$, and with $h \rightarrow 0$ we can choose $\Omega'$ so that $\lambda^d(\Omega') \rightarrow \lambda^d(\Omega)$. Moreover, $|\nabla I_h v | < C$ almost everywhere, and hence discontinouities do not cause problems. It follows that $\norm{\nabla (I_h v - v)}_{L^2} \rightarrow 0$.
\end{proof}

\end{appendix}


\begin{thebibliography}{10}

\bibitem{bai2000templates}
Zhaojun Bai, James Demmel, Jack Dongarra, Axel Ruhe, and Henk van~der Vorst.
\newblock {\em Templates for the solution of algebraic eigenvalue problems: a
  practical guide}.
\newblock SIAM, 2000.

\bibitem{quadraturepaper}
Uday Banerjee and John~E Osborn.
\newblock Estimation of the effect of numerical integration in finite element
  eigenvalue approximation.
\newblock {\em Numerische Mathematik}, 56(8):735--762, 1989.

\bibitem{Bogacki1996}
P.~Bogacki and L.~F. Shampine.
\newblock An efficient {R}unge--{K}utta (4, 5) pair.
\newblock {\em Computers \& Mathematics with Applications}, 32(6):15--28, 1996.

\bibitem{cai95}
Xiao-Chuan Cai.
\newblock The use of pointwise interpolation in domain decomposition methods
  with nonnested meshes.
\newblock {\em SIAM Journal on Scientific Computing}, 16(1):250--256, 1995.

\bibitem{juafem}
K.~Carlsson.
\newblock {KristofferC/JuAFEM.jl}: {F}inite element toolbox for {J}ulia, 2019.

\bibitem{davies1982pt2}
Edward~Brian Davies.
\newblock Metastable states of symmetric {M}arkov semigroups {II}.
\newblock {\em Journal of the London Mathematical Society}, 2(3):541--556,
  1982.

\bibitem{DeJu99}
Michael Dellnitz and Oliver Junge.
\newblock On the approximation of complicated dynamical behavior.
\newblock {\em SIAM Journal on Numerical Analysis}, 36(2):491--515, 1999.

\bibitem{deuflhard2005robust}
Peter Deuflhard and Marcus Weber.
\newblock Robust {P}erron cluster analysis in conformation dynamics.
\newblock {\em Linear algebra and its applications}, 398:161--184, 2005.

\bibitem{ernguermond}
Alexandre Ern and Jean-Luc Guermond.
\newblock {\em Theory and practice of finite elements}, volume 159.
\newblock Springer Science \& Business Media, 2013.

\bibitem{evans}
L.C. Evans.
\newblock {\em Partial Differential Equations}.
\newblock Graduate studies in mathematics. American Mathematical Society, 1997.

\bibitem{froyland05}
Gary Froyland.
\newblock Statistically optimal almost-invariant sets.
\newblock {\em Physica D: Nonlinear Phenomena}, 200(3-4):205--219, 2005.

\bibitem{froyland15}
Gary Froyland.
\newblock {Dynamic isoperimetry and the geometry of Lagrangian coherent
  structures}.
\newblock {\em Nonlinearity}, 28(10):3587, 2015.

\bibitem{rbf}
Gary Froyland and Oliver Junge.
\newblock On fast computation of finite-time coherent sets using radial basis
  functions.
\newblock {\em Chaos: An Interdisciplinary Journal of Nonlinear Science},
  25(8):087409, 2015.

\bibitem{froylandjunge18}
Gary Froyland and Oliver Junge.
\newblock Robust {FEM}-based extraction of finite-time coherent sets using
  scattered, sparse, and incomplete trajectories.
\newblock {\em SIAM Journal on Applied Dynamical Systems}, 17(2):1891--1924,
  2018.

\bibitem{FK17}
Gary Froyland and Eric Kwok.
\newblock A dynamic {L}aplacian for identifying {L}agrangian coherent
  structures on weighted {R}iemannian manifolds.
\newblock {\em Journal of Nonlinear Science}, page To appear, 2017.

\bibitem{froyland2010coherent}
Gary Froyland, Simon Lloyd, and Naratip Santitissadeekorn.
\newblock Coherent sets for nonautonomous dynamical systems.
\newblock {\em Physica D: Nonlinear Phenomena}, 239(16):1527--1541, 2010.

\bibitem{froylandpadberg09}
Gary Froyland and Kathrin Padberg.
\newblock Almost-invariant sets and invariant manifolds—connecting
  probabilistic and geometric descriptions of coherent structures in flows.
\newblock {\em Physica D: Nonlinear Phenomena}, 238(16):1507--1523, 2009.

\bibitem{froyland2019sparse}
Gary Froyland, Christopher~P Rock, and Konstantinos Sakellariou.
\newblock Sparse eigenbasis approximation: Multiple feature extraction across
  spatiotemporal scales with application to coherent set identification.
\newblock {\em Communications in Nonlinear Science and Numerical Simulation},
  77:81--107, 2019.

\bibitem{karraschkeller}
Daniel Karrasch and Johannes Keller.
\newblock {A geometric heat-flow theory of Lagrangian coherent structures}.
\newblock {\em arXiv preprint arXiv:1608.05598}, 2017.

\bibitem{kato}
Tosio Kato.
\newblock {\em Perturbation theory for linear operators}, volume 132.
\newblock Springer Science \& Business Media, 2013.

\bibitem{voronoidelaunay}
A.~Keselman et~al.
\newblock {JuliaGeometry/VoronoiDelaunay.jl}: {F}ast and robust {V}oronoi \&
  {D}elaunay tesselation creation with {J}ulia, 2019.

\bibitem{juliadiffeq}
C.~Rackauckas and Q.~Nie.
\newblock {DifferentialEquations.jl} -- a performant and feature-rich ecosystem
  for solving differential equations in julia.
\newblock {\em Journal of Open Research Software}, 5(1):15, 2017.

\bibitem{rypina2007lagrangian}
II~Rypina, Michael~G Brown, Francisco~J Beron-Vera, Huseyin Ko{\c{c}}ak,
  Maria~J Olascoaga, and IA~Udovydchenkov.
\newblock {On the Lagrangian dynamics of atmospheric zonal jets and the
  permeability of the stratospheric polar vortex}.
\newblock {\em Journal of the Atmospheric Sciences}, 64(10):3595--3610, 2007.

\bibitem{strangfix}
Gilbert Strang and George~J Fix.
\newblock {\em An analysis of the finite element method}, volume 212.
\newblock Prentice-Hall Englewood Cliffs, NJ, 1973.

\end{thebibliography}
\end{document}